\documentclass[11pt,a4paper,reqno]{amsart}
\usepackage{version}
\usepackage{fancyhdr}
\usepackage{color}   
\usepackage{hyperref}
\hypersetup{
    colorlinks=true, 
    linktoc=all,     
    linkcolor=blue,  
    citecolor=red,
}

\usepackage{romannum}
\AtBeginDocument{\pagenumbering{arabic}}

\usepackage{xcolor}

\usepackage{amsmath,amsthm,amssymb,latexsym,epsfig,graphicx,subfigure}
\usepackage{color}
\usepackage{parskip}
\usepackage{hyperref}
\usepackage{mathabx}
\numberwithin{equation}{section}

\newtheorem{theorem}{Theorem}[section]

\newtheorem{lemma}{Lemma}[section]
\newtheorem{prop}{Proposition}[section]

\newtheorem{corollary}{Corollary}[section]
 
\theoremstyle{definition}
\newtheorem{definition}[theorem]{Definition}
\theoremstyle{definition}
\newtheorem{remark}{Remark}[section]

\newcommand{\R}{\mathbb{R}}

\newcommand{\Z}{\mathbb{Z}}
\newcommand{\N}{\mathbb{N}}

\newcommand{\cI}{\mathcal{I}}

\newcommand{\e}{\varepsilon}

\newcommand{\diam}{\textrm{diam}}

\newcommand{\spt}{\textrm{spt }}

\usepackage{tikz}
\usepackage{graphicx}

\usepackage{comment}

\newcommand{\norm}[1]{{\left\Vert#1\right\Vert}}

\allowdisplaybreaks

\begin{document}
\title{\texorpdfstring{$L^{p}$}{Lp}-integrability of functions with Fourier supports on fractal sets on the moment curve}
\author{Shengze Duan}
\address{Department of Mathematics, University of Rochester, USA.} \email{sduan4@ur.rochester.edu}
\author{Minh-Quy Pham}
\address{Department of Mathematics, University of Rochester, USA.} \email{qpham3@ur.rochester.edu}
\author{Donggeun Ryou}
\address{Department of Mathematics, Indiana University Bloomington, USA.} \email{dryou@iu.edu }

\begin{abstract} For $0 < \alpha \leq 1$, let $E$ be a compact subset of the $d$-dimensional moment curve in $\R^d$ such that $N(E,\e) \lesssim \e^{-\alpha}$ for $0 <\e <1$ where $N(E,\e)$ is the smallest number of $\e$-balls needed to cover $E$. We proved that if $f \in L^p(\R^d)$ with
\begin{align*}
        1 \leq p\leq p_\alpha:=
        \begin{cases}
            \frac{d^2+d+2\alpha}{2\alpha} & d \geq 3,\\
            \frac{4}{\alpha} &d =2,
        \end{cases}
    \end{align*}
and $\widehat{f}$ is supported on the set $E$, then $f$ is identically zero. We also proved that the range of $p$ is optimal by considering random Cantor sets on the moment curve. We extended the result of Guo, Iosevich, Zhang and Zorin-Kranich \cite{Guoetal23}, including the endpoint. We also considered applications of our results to the failure of the restriction estimates and Wiener Tauberian Theorem.
\end{abstract}
\maketitle
\tableofcontents

\section{Introduction}

This paper concerns the following question:

Assume that $f \in L^p(\R^d)$ and its Fourier transform is supported on the set $E$ in $\R^d$. What is the optimal range of the $p$ such that $f$ is identically zero?  

Agranovsky and Narayanan \cite{Agranovsky04} proved that if $E$ is a $C^1$ manifold of dimension $0 < \alpha < d$ and $1 \leq p \leq 2d/\alpha$, then $f $ is identically zero. Their result is optimal if $\alpha \geq d/2$. A typical example is when the manifold is the $d-1$ dimensional paraboloid. If $\mu$ is the surface measure on the paraboloid, we let $f = \widehat{\mu}$. Then $f$ is nonzero and $f \in L^p(\R^d)$ for $p > 2d/(d-1)$. Their result was extended to fractal sets by Senthil Raani \cite{Raani14}. If the $\alpha$-dimensional packing pre-measure of $E$ is finite, the optimal range of $p$ where $f\equiv 0$ is $1 \leq p \leq 2d/\alpha$ (See Remark \ref{premeas}).

If we assume that the $E$ is a compact set with finite $\alpha$-dimensional Hausdorff measure, then the optimal range of $p$ such that $f \equiv 0$ is $1 \leq p < 2d/\alpha$. Thus, it does not include the endpoint $2d/\alpha$. Salem \cite{Salem51} proved the sufficiency of the range of $p$ when $d=1$, and Edgar and Rosenblatt \cite{ER79} proved it when $\alpha \geq d-1$. For $0 < \alpha <d$, the proof of the sufficiency of the range of $p$ follows from the work of Kahane \cite{Kahane84} (see \cite{Dob24}). The necessity of the range of $p$ was recently proved by Dobronravov \cite{Dob24}.

However, the threshold $2d/\alpha$ is not always optimal. The moment curve 
$$\Gamma_d=\{\gamma_d(t)=(t,t^2,\dots, t^d): t\in [0,1]\}$$ 
is such an example. Recently, Guo, Iosevich, Zhang and Zorin-Kranich \cite{Guoetal23} proved that if $E = \Gamma_d$ and if 
$$1 \leq p < \frac{d^2+d+2}{2},$$
then $f \equiv 0$. The range of $p$ is optimal except for the endpoint due to the work of Arkhipov, Chubarikov, and Karatsuba \cite[Theorem 1.3]{ACK_book}.

The main aims of this paper are the following. First, we extend the result of \cite{Guoetal23} to fractal sets on $\Gamma_d$ including the endpoint. Second, we apply our result to the necessity condition of restriction estimates on the moment curve and Wiener Tauberian Theorem.

Throughout this paper, we denote by $A \lesssim B$ when $A \leq CB$ for some constant $C>0$, and we denote by $A \sim B$ when $A \lesssim B$ and $B \lesssim A$. If the constant $C$ depends on a parameter such as $\e$, we write $A(\e) \lesssim_\e B(\e)$. The Lebesgue measure of a set $A$ in $\R^d$ is denoted by $|A|$, and the cardinality of a finite set $B$ is denoted by $\#B$. For $x \in \mathbb{R}$, we denote by $e(x)$ the quantity $e^{2\pi i x}$.

\subsection{\texorpdfstring{$L^p$}{Lp}-integrability of functions with Fourier supports on fractal sets on the moment curve}
Let $E\subset \R^d$ be a non-empty bounded set. For $0<\e<1$, let {\it$\e$-covering number} $N(E,\e)$ be the smallest number of $\e$-balls needed to cover $E$:
\begin{align*}
    N(E,\e)=\min\{ k: E\subset \bigcup\limits_{i=1}^kB(x_i,\e), \text{ for $x_i\in \R^d$}\}.
\end{align*}

Let us state our first main result.
\begin{theorem}\label{thm_main_Lp_moment}
    Let $d\geq 2$ and $0< \alpha\leq 1$. Let $E$ be a set on $\Gamma_d$ such that $N(E,\e)\lesssim \e^{-\alpha}$ uniformly in $0 < \e <1$. If $f\in L^p(\R^d)$ and $\spt \widehat{f}\subseteq E$, then $f \equiv 0$ when
    \begin{align}\label{eq_thm_1}
        1\leq p\leq p_\alpha:=
        \begin{cases}
            \frac{d^2+d+2\alpha}{2\alpha} & \text{if} \ d \geq 3,\\
            \frac{4}{\alpha} &\text{if} \ d =2.
        \end{cases}
    \end{align}
\end{theorem}

The range of $p$ in Theorem \ref{thm_main_Lp_moment} is optimal in the following sense.
\begin{theorem}\label{thm_main_Lp_optimal}
    Let $d\geq 2$, and $0< \alpha\leq 1$. If $p  > p_\alpha$, there exists a nonzero $f \in L^p(\R^d)$ such that  $N(\spt {\widehat{f}}, \e) \lesssim \e^{-\alpha} $ for all $ 0< \e <1$. 
\end{theorem}

\begin{remark}\label{premeas}
    Theorem 2.3 in \cite{Raani14} was stated with an assumption on the $\alpha$-dimensional packing measure. However, the author clarified that it should be the $\alpha$-dimensional packing pre-measure \cite{Raani_pri}. For $\alpha \geq 0$ and $\delta >0$, let
    \[
    \mathcal{P}_\delta^\alpha(E) = \sup\ \sum_{i=1}^\infty |B_i|^\alpha,
    \]
    where the supremum is taken over all collection of disjoint balls $\{B_i\}$ of raddi at most $\delta$ with centers in $E$. The $\alpha$-dimensional packing pre-measure $\mathcal{P}_0^\alpha(E)$ is defined by 
    \[
    \mathcal{P}_0^\alpha(E) = \lim_{\delta \rightarrow 0} \mathcal{P}_\delta^{\alpha}(E). 
    \]
    Then, Lemma 2.2 in \cite{Raani14} holds when $\mathcal{P}_0^\alpha(E) < \infty $ and so does Theorem 2.3. 
\end{remark}
\begin{remark}\label{re_covnum}
    Lemma 2.2 in \cite{Raani14} is still valid with a slightly weaker condition that $N(E,\e)  \lesssim \e^{-\alpha}$ for all $0 < \e < 1$ instead of $\mathcal{P}_0^\alpha(E) <\infty$, since $E(\e) \lesssim \e^{n} N(E,\e) \lesssim \e^{n-\alpha}$ where $E(\e ) = \{x \in \R^d: d(x,E) < \e \} $.
    Therefore, Theorem 2.3 in \cite{Raani14} can be restated as follows: 
    
    Let $f \in L^p(\R^d)$ be a function such that $\widehat{f}$ is supported on a set $E$ with $N(E,\e) \lesssim \e^{-\alpha}$ for all $0< \e <1$. Then $f \equiv 0$, provided that $p \leq 2d/\alpha$.

    We will use this when $d=2$.
\end{remark}

\begin{remark}
    In fact, the proof of Theorem \ref{thm_main_Lp_moment} for $d \geq 3$ works even when $d=2$ and Senthil Raani's result can be applied to higher dimensions. By choosing a better one between the two, $p_\alpha$ can be also written as
    \[
    p_\alpha = \max\left\{ \frac{2d}{\alpha} , \frac{d^2+d+2\alpha}{2\alpha} \right\}.
    \]
    If $d=2$, then $\frac{2d}{\alpha}$ is larger and if $d \geq 3$, then $\frac{d^2+d+2\alpha}{2\alpha}$ is larger.
\end{remark}

\subsection{Restriction estimates on the moment curve}
In this section, we consider the extension estimate
\begin{equation}\label{restriction1}
\norm{\widehat{fd\mu}}_{L^p(\R^d)} \lesssim \norm{f}_{L^q(\mu)}    
\end{equation}
where $\mu$ is a measure supported on $\Gamma_d$.

Since its dual version is called the restriction estimate, we will describe our results by using the extension estimate. If $\mu$ is the arc length measure on $\Gamma_d$, i.e.,
\[
\int f d\mu = \int f(t, t^2, \cdots, t^d) dt,
\]
then the optimal range of $p$ and $q$ where the extension estimate holds is already known.
\begin{theorem}\label{restric_thm_moment}
    Let $\mu$ be the arc length measure on $\Gamma_d$ defined above. Then \eqref{restriction1} holds if and only if 
    \[
p \geq q'\frac{d(d+1)}{2} \qquad \text{and} \qquad p > \frac{d^2+d+1}{2},
\]
where $q'$ is the conjugate of $q$, i.e. $1/q+1/q'=1$.
\end{theorem}
When $d=2$ and $d=3$, the result was proved by Zygmund \cite{Zygmund} and Prestini \cite{Prestini} respectively. In higher dimensions, it was shown by Drury \cite{Drury}. 

The necessity of the first restriction is due to the Knapp-type example. For example, see Example 1.8 in \cite{Dem_book}. The necessity of the second restriction follows by using a constant function on the support of $\mu$, which can be found in Theorem 1.3 in \cite{ACK_book}. For more general discussion, see \cite{BGGIST07}.

As an analog of the necessity of the range of $p$ and $q$ in Theorem \ref{restric_thm_moment}, we have the following.
\begin{theorem}\label{restriction_cant_any}
    Assume that $\mu$ is a nonzero Borel probability measure with $\spt \mu \subseteq \Gamma_d$. If $N(\spt \mu ,\e) \lesssim \e^{-\alpha}$ for all $0 < \e <1$, then \eqref{restriction1} cannot hold when
    \begin{equation}\label{necessity_2}
        1 \leq p \leq p_\alpha.
    \end{equation}
    If $\mu(B(x,r)) \lesssim r^\alpha$ for all $x \in \R^d$ and $r>0$, then \eqref{restriction1} cannot hold when
    \begin{equation}\label{necessity_1}
    p < q' \frac{d(d+1)}{2\alpha}.    
    \end{equation}
    
\end{theorem}
Note that if the measure $\mu$ is AD-regular, i.e. $\mu(B(x,r)) \sim r^\alpha$ for all $x \in \spt\mu$ and $0<r<\diam(\spt \mu)$, then the measure also satisfies that $N(\spt \mu,\e) \lesssim \e^{-\alpha}$ for all $ 0 < \e <1$.

The necessity of the first restriction \eqref{necessity_2} follows from a stronger result.
\begin{prop} 
    Let $\mu$ be a nonzero Borel measure with $\spt \mu\subseteq \Gamma_d$. If $N(\spt \mu,\e)\lesssim \e^{-\alpha}$ for all $ 0 <\e <1 $ and $p \leq p_\alpha$, then \eqref{restriction1} cannot be satisfied for any nonzero $f \in L^q(\mu)$. 
\end{prop}
\begin{proof}
    If there exists such $f \in L^q(\mu)$, then $\widehat{f{d\mu} } (-x) \in L^p(\R^d)$.
    This contradicts Theorem \ref{thm_main_Lp_moment} if we replace $f$ in Theorem \ref{thm_main_Lp_moment} by $\widehat{fd\mu}(-x)$.

\end{proof}

Thus, it suffices to prove the necessity of the second restriction \eqref{necessity_1} and it can be shown by a Knapp-type example.
\begin{definition}\label{dualbox}
    Let $B$ be a rectangular box in $\R^d$ with side lengths $l_1, \cdots, l_d$. The \textit{dual rectangular box $\widetilde{R}$} of the box $R$ is defined by the rectangular box with side lengths $l_1^{-1}, \cdots, l_d^{-1} $ centered at the origin such that $l_1$ side of $R$ is parallel to the $l_1^{-1}$ side of $\widetilde{R}$.
\end{definition}
\begin{proof}[Proof of Theorem \ref{restriction_cant_any}]
    The argument in Proposition 3.1 in \cite{Mitsis02} shows that if $\mu(B,r) \lesssim r^\alpha$ for all $r>0$, then there exist sequences $\{x_k\}_{k \in \N}$ and $\{r_k\}_{k \in \N}$ such that $r_k \rightarrow 0$ as $k \rightarrow \infty$ and 
    $$
    \mu(B(x_k,r_k)) \sim r_k^\alpha.
    $$
   For each $k$, let be $f_k$ the indicator function of $B(x_k,r_k)$. Since $\mu$ is supported on $\Gamma_d$, $B(x_k,r_k) \cap \Gamma_d$ is contained in a rectangular box $B_k$ with sides $\sim r_k, \cdots, r_k^{d}$. Then, let $\widetilde{B}_k$ be its dual rectangular box. Note that $\widetilde{B}_k
   $ has side lengths $r_k^{-1}$, $r_k^{-2}, \cdots, r_k^{-d}$. Since $B(x_k, r_k) \cap \Gamma_d \subseteq B_k$, if $\xi \in \widetilde{B}_k/100$, then $e(-\xi\cdot x) \sim 1 $, i.e., $\text{Re} (e(-\xi\cdot x)) \sim 1 $ and $\text{Im}(e(-\xi\cdot x)) \leq 1/10 $. Therefore, we obtain that
    \begin{align*}
    \norm{\widehat{f_kd\mu}}_{L^p(\R^d)}^p&\geq \int_{\widetilde{B}_k/100} \left|\int_{B(x_k,r_k)}e(-\xi \cdot x) d\mu(x)\right|^pd\xi\\
    &\gtrsim \int_{\widetilde{B}_k/100} \left| \mu(B(x_k, r_k))\right|^pd\xi\\
    &\gtrsim r_k^{\alpha p - \frac{d(d+1)}{2}}.
    \end{align*}
Also, we have that
$$ \norm{f_k}_{L^q(\mu)} \lesssim r_k^{\frac{\alpha}{q}}.$$

If we assume \eqref{restriction1}, we combine the estimates above and obtain that
\[
r_k^{\alpha -\frac{d(d+1)}{2p}} \lesssim r_k^{\frac{\alpha}{q}}.
\]
Since $r_k \rightarrow 0$ as $k \rightarrow \infty$, we have the desired relation between $p$ and $q$. 
\end{proof}

\begin{figure}[ht]
    \centering

\tikzset{every picture/.style={line width=0.75pt}} 

\begin{tikzpicture}[x=0.75pt,y=0.75pt,yscale=-1,xscale=1]

\draw  [draw=none][fill={rgb, 255:red, 155; green, 155; blue, 155 }  ,fill opacity=1 ] (331.5,289.56) -- (181,217.56) -- (331.5,217.56) -- cycle ;
\draw    (109,301.56) -- (109,11.06) ;
\draw [shift={(109,9.06)}, rotate = 90] [color={rgb, 255:red, 0; green, 0; blue, 0 }  ][line width=0.75]    (10.93,-3.29) .. controls (6.95,-1.4) and (3.31,-0.3) .. (0,0) .. controls (3.31,0.3) and (6.95,1.4) .. (10.93,3.29)   ;
\draw    (97,289.56) -- (386.5,289.07) ;
\draw [shift={(388.5,289.06)}, rotate = 179.9] [color={rgb, 255:red, 0; green, 0; blue, 0 }  ][line width=0.75]    (10.93,-3.29) .. controls (6.95,-1.4) and (3.31,-0.3) .. (0,0) .. controls (3.31,0.3) and (6.95,1.4) .. (10.93,3.29)   ;
\draw  [fill={rgb, 255:red, 155; green, 155; blue, 155 }  ,fill opacity=1 ] (109,67.06) -- (331.5,67.06) -- (331.5,218.06) -- (109,218.06) -- cycle ;
\draw    (331.5,218.06) -- (331.5,289.56) ;
\draw  [dash pattern={on 4.5pt off 4.5pt}]  (181,218.06) -- (181,289.56) ;
\draw  [dash pattern={on 4.5pt off 4.5pt}]  (331.5,289.56) -- (181,217.56)  ;

\draw (84.5,4.4) node [anchor=north west][inner sep=0.75pt]    {$\frac{1}{p}$};
\draw (406,266.9) node [anchor=north west][inner sep=0.75pt]    {$\frac{1}{q}$};
\draw (325.5,300.9) node [anchor=north west][inner sep=0.75pt]    {$1$};
\draw (90,59.4) node [anchor=north west][inner sep=0.75pt]    {$1$};
\draw (81,195.4) node [anchor=north west][inner sep=0.75pt]    {$\frac{1}{p_{\alpha }}$};
\draw (170,299.4) node [anchor=north west][inner sep=0.75pt]    {$\frac{1}{q_{\alpha }}$};
\draw (210,136) node [anchor=north west][inner sep=0.75pt]   [align=left] {$A$};
\draw (281.5,231.5) node [anchor=north west][inner sep=0.75pt]   [align=left] {$B$};

\end{tikzpicture}

    \caption{Failure range of the extension estimate: $q_\alpha=
        \begin{cases}
            p_\alpha & d \geq 3\\
            4 &d =2.
        \end{cases}$}
    \label{fig1}
\end{figure}
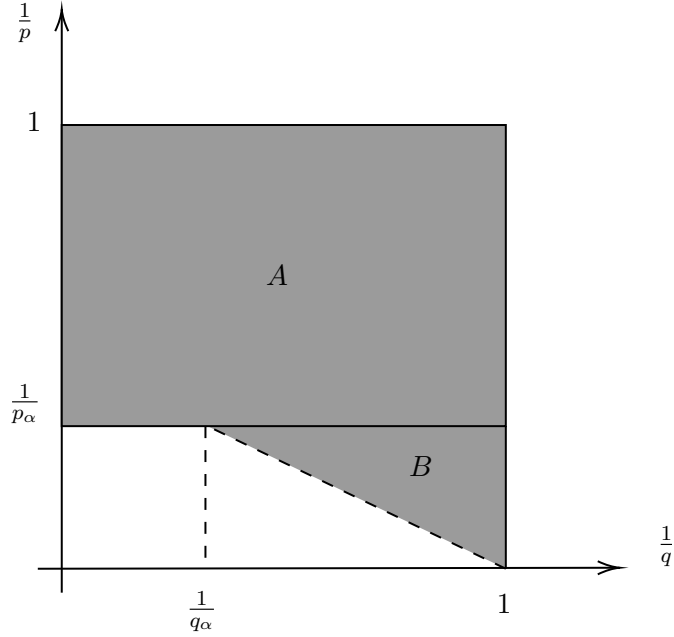

Therefore, if the measure $\mu$ is AD-regular, Theorem \ref{restriction_cant_any} implies that we can divide the range of $(p,q)$ where \eqref{restriction1} fails, into two parts; when \eqref{necessity_2} holds and when \eqref{necessity_1} holds but \eqref{necessity_2} does not hold. They correspond to the region $A$ and the region $B$ respectively in Figure \ref{fig1}. In the region $A$, the estimate \eqref{restriction1} fails for any nonzero $f \in L^q(\mu)$ and in the region $B$, the estimate \eqref{restriction1} fails for some $f \in L^q(\mu)$. It would be interesting to find for which $f \in L^q(\mu)$ the estimate \eqref{restriction1} holds or fails in the region $B$.

Also, when $\alpha =1$, the arc length measure on $\Gamma_d$ is an example such that \eqref{restriction1} holds if and only if 
\begin{equation}\label{sufficient_1}
p \geq q'\frac{d(d+1)}{2 \alpha} \qquad \text{and} \qquad p > p_\alpha.    
\end{equation}
However, when $0< \alpha <1$, there is no such example known yet. When $d=2$, Ryou \cite{Ryou23} constructed a measure such that \eqref{restriction1} holds when $p > 6/\alpha$ and $q=2$. Even with interpolation, this result only covers a partial range of \eqref{sufficient_1}. Thus, it would be interesting to construct a measure $\mu$ on $\Gamma_d$ with $0 < \alpha <1$ such that \eqref{restriction1} holds if and only if $p$ and $q$ satisfy \eqref{sufficient_1}.

\subsection{Wiener Tauberian Theorem}
Consider a function $f \in L^p(\R^d)$. Let $\mathcal{M}_f$ denote the space of the finite linear combinations of translates of $f$. Wiener Tauberian Theorem concerns a necessary and sufficient condition of $f$ such that $\mathcal{M}_f$ is dense in $L^p(\R^d)$. Wiener \cite{Weiner} proved the following result. See also p. 234 in  \cite{Dono}.
\begin{theorem}
    If $f \in L^1(\R^d)$, then $\mathcal{M}_f$ is dense in $L^1(\R^d)$ if and only if the zero set of $\widehat{f}$ is empty.\\
    If $f \in L^2(\R^d)$, then $\mathcal{M}_f$ is dense in $L^2(\R^d)$ if and only if the zero set of $\widehat{f}$ has Lebesgue measure zero.
\end{theorem}
Wiener conjectured that a similar result would hold when $1 < p<2$. However, it was disproved by Lev and Olevskii \cite{LevOl}. They showed that, if $1 <  p < 2$, then there exists two functions $f$ and $g$ in $L^1(\R) \cap C_0 (\R)$ such that $\mathcal{M}_f$ is dense in $L^p(\R)$ but $\mathcal{M}_g $ is not, while the zero sets of $\widehat{f}$ and $\widehat{g}$ are the same.

However, there are some results about the sufficient conditions of $f$. Beurling \cite{Beur} proved that if $f \in L^p(\R)$ and the zero set of $\widehat{f}$ is in the closed set of Hausdorff dimension $0< \alpha< 1 $, then $\mathcal{M}_f$ is dense in $L^p(\R)$ for $ 2/(2-\alpha) <p <\infty$. 

Herz \cite{Herz} and Agranovsky and Narayanan  \cite{Agranovsky04} proved similar results in higher dimensions with additional assumptions. Senthil Raani \cite{Raani14} improved their work with an additional hypothesis on the zero set of $\widehat{f}$ as follows.
\begin{prop}\label{Winer_Raani}
    Assume that $f \in L^{p'} \cap L^1(\R^d) $ and $1 \leq p < \infty$ where $1/p +1/p'=1$. If the zero set of $\widehat{f}$ is contained in $E \subseteq \Gamma_d$ with $p \leq 2d/\alpha$ and $N(E,\e) \lesssim \e^{-\alpha}$ for $0<\e< 1$, then $\mathcal{M}_f$ is dense in $L^{p'}(\R^d)$.
\end{prop} 
Using Theorem \ref{thm_main_Lp_moment}, we can prove an analog of the above result.

\begin{prop}\label{Wiener_moment}
    Assume that $f\in L^{p'}\cap L^1$. If the zero set of $\hat{f}$ is contained in $ E\subseteq\Gamma_d$ with $p\leq p_{\alpha}$ and $N(E,\e)\lesssim \e^{-\alpha}$ for all $0 < \e <1$, then $\mathcal{M}_f$ is dense in $L^{p'}(\R^d)$.
\end{prop}

The proof of Proposition \ref{Wiener_moment} is the same as Proposition \ref{Winer_Raani} except that we use Theorem \ref{thm_main_Lp_moment}. Since $p_\alpha > 2d/\alpha $ for $d \geq 3$, the range of $p'$ in Proposition \ref{Wiener_moment} is larger than the range of $p'$ in Proposition \ref{Winer_Raani}.

\subsection{Outline of the Paper} The rest of the paper is devoted to the proof of Theorems \ref{thm_main_Lp_moment} and \ref{thm_main_Lp_optimal}. In Section \ref{sec2}, we prove Theorem \ref{thm_main_Lp_moment}. We follow the argument in the paper of Senthil Raani \cite{Raani14} with some modifications. We use rectangular boxes with side lengths $\e, \e^{2}, \cdots, \e^d$ adapted to the moment curve instead of a ball of radius $\e$. In Section \ref{sec3}, we prove Theorem \ref{thm_main_Lp_optimal} by showing the $L^p$ integrability of a random Cantor set on the moment curve. 

To the authors' best knowledge, there is no explicit example that achieves Theorem 1.2 except for $\alpha=1$. If $\alpha<1$, we need to consider fractal measures, and classical theories on oscillatory integrals cannot be directly applied. Thus, we consider a random Cantor set, which allows us to combine it with Hoeffding's inequality.

The Fourier decay of a random Cantor set was studied by Shmerkin and Suomala \cite{ShSu17} (also, see \cite{Shsu18}). Ryou \cite{Ryou23} combined the idea from \cite{ShSu17} and the estimates on the oscillatory integral, and then found an estimate on the Fourier decay of a random Cantor set constructed on the parabola. We will use this result when $d=2$. 

If $d \geq 3$ and $\mu$ is the arc length measure on $\Gamma_d$, it is known that $|\widehat{\mu}(\xi)| \lesssim |\xi|^{-1/d}$. However, $\widehat{\mu}(\xi)$ decays faster than $|\xi|^{-1/d}$ in most directions of $\xi$. More precisely, a subset of directions of $\xi$ in the $d-1$ dimensional sphere $S^{d-1}$ such that $|\widehat{\mu}(\xi)| \sim |\xi|^{-1/d} $ has dimension $1$. Such a subset is significantly small when $d \geq 3$. Therefore, we consider a decomposition according to how far $\xi$ is from the `bad' directions, where the slowest Fourier decay occurs.\\
Another difficulty is that fractal measures are supported on many small, disjoint intervals. If $N(E, \e) \lesssim \e^{-\alpha}$, we need to consider the summation of $\sim \e^{-\alpha}$ many oscillatory integrals on the disjoint intervals of length $\e$. The same issue was also considered in \cite{Ryou23}. When $d\geq 3$, we use the result of Arkhipov, Chubarikov, and Karatsuba \cite{ACK_book} on each oscillatory integral. Then, by using the ideas in \cite{Ryou23} and \cite{ShSu17}, we show that the sum of these oscillatory integrals satisfies the desired upper bound on each decomposition.

\subsection{Acknowledgements}
We would like to thank Alex Iosevich for many discussions and hence, for inspiring this work. We would also like to thank Allan Greenleaf and K.S. Senthil Raani for helpful conversations.

\section{Proof of Theorem \ref{thm_main_Lp_moment}}\label{sec2}
\
In this section, we will give the proof of Theorem \ref{thm_main_Lp_moment}. The proof is a modification of the arguments in \cite{Guoetal23} and \cite{Raani14}. 

For $t\in \R$, let $(\Vec{e_1}(t),\dots, \Vec{e_d}(t))$ denote the Frenet coordinates along the moment curve at the point $\gamma_d(t)$. For $\e>0$, define an $\e$-isotropic neighborhood of the moment curve at the point $\gamma_d(t)$ by
    \begin{align}\label{def_n_gamma}
        \Gamma_{\e,t}:=\{\gamma_d(t)+\e_1\Vec{e_1}(t)+\cdots+\e_d\Vec{e_d}(t): |\e_k|\leq \e^k, 1\leq k\leq d\}.
    \end{align}
Observe that by the assumption, for $\e>0$ small enough, we then can cover $ E$ by a finite-overlapping family of rectangular boxes $\{B_{\e,i}\}_{i=1}^{M_\e}$ with $M_\e \lesssim \e^{-\alpha}$ such that $B_{\e,i} = \Gamma_{\e,t_i}$ and the points $t_i$ are $\sim \e$-separated. Note that each $B_{\e,i}$ has side lengths $\e, \e^2, \cdots, \e^d$.

    Let $\{\phi_i^{(1)}\}_{i=1}^{M_\e}$ be a partition of unity in the first variable $x_1$ such that each $\phi_i^{(1)}$ is supported on an interval of length $\sim \e$ and the center of $\phi_i^{(1)}$ is the same with the first coordinate of $B_{\e,i}$ respectively. Then, we choose a smooth function $\phi^{(2)}(x')$ where $x' \in \R^{d-1}$ such that $\phi^{(2)}(x') =1$ for $x' \in B(0,2)$ and let $\phi_i(x) = \phi_i^{(1)}(x_1)\phi^{(2)}(x')$.
    
    For each $i$, we define
    \begin{align*}
        f_{\e,i}=\mathcal{F}^{-1}\big(\widehat{f}\phi_i\big).
    \end{align*}
    where $\mathcal{F}^{-1}$ is the inverse Fourier transform. In other words, $\widehat{f_{\e,i}}=\widehat{f}\phi_{i}$. We multiplied the smooth function $\phi_i$, since $\widehat{f}$ is well defined as a tempered distribution. Let $CB_{\e,i}$ be a rectangle with side lengths $C\e, C\e^2, \cdots, C\e^d$ and the same center as $B_{\e,i}$.
    Note that $\widehat{f_{\e,i}}$ is supported on $C B_{\e,i}$ for some constant $C$ independent on $i $ and $\e$ and also we have $f = \sum_{i}f_{\e,i}$.

    Let $\phi$ be a normalized smooth function supported on $B(0,1)$. We also define that $D_\e$ is the diagonal matrix with diagonal entries $\e, \cdots, \e^d$ and $M_t$ is the matrix whose columns are 
    $\Vec{e_1}(t),\cdots,\Vec{e_d}(t)$. For each $i$, define 
    \begin{align*}
        \eta_{\e,i}(x)=\e^{-\frac{d^2+d}{2}}\phi(D_\e^{-1} M_{t_i}^{-1} x).
    \end{align*}
    
    Denoting $u=\widehat{f}$, we then define
    \begin{align*}
        u_\e=\sum\limits_{i=1}^{M_\e} \widehat{f_{\e,i}}\ast \eta_{\e,i}=\sum\limits_{i=1}^{M_\e} \left(u{\phi_{i}}\right)\ast \eta_{\e,i}.
    \end{align*}
    To prove the theorem, we need the following lemmas. 
    \begin{lemma}\label{lemma_1}
        Assume that $f \in L^{p_\alpha}(\R^d)$. For $0 < \e <1$, there exist $\{a_j\}_{j\in \Z}$ and $\{b_{j,\e}\}_{j\in \Z}$ such that
        \begin{align*}
        \Vert u_\e\Vert_2^2 \lesssim \e^{-\frac{d^2+d-2\alpha}{2}}\sum\limits_{j=-\infty}^{\infty} a_jb_{j,\e},
    \end{align*}
    where $\sum_j|a_j|<\infty$, and for any fixed $j\in\Z$, $|b_{j,\e} | \lesssim \norm{f}_{p_\alpha}^2$ and $b_{j,\e}\to 0$, as $\e\to 0$.
    \end{lemma}
    \begin{lemma}\label{lemma_2}
    Assume that $f \in L^p(\R^d)$ for some $p \geq 2$. For any compactly supported smooth function $\psi:\R^d\to \R$, we have
    \begin{align*}
        |\langle u, \psi\rangle|^2
    &=\lim\limits_{\e\to 0}|\langle u_\e, \psi\rangle|^2.
    \end{align*}
\end{lemma}

Assuming the above lemmas, let us prove Theorem \ref{thm_main_Lp_moment} first.
\begin{proof}[Proof of Theorem \ref{thm_main_Lp_moment} using Lemmas \ref{lemma_1} and \ref{lemma_2}]  
As it was pointed out in Remark \ref{re_covnum}, the case when $d=2$ is a direct result of Senthil Raani \cite{Raani14}. Thus, we consider only when $d \geq 3$. Also, we will only consider when $p = p_\alpha$. If $f \in L^p(\R^d)$ for $p < p_\alpha$, we can convolve $f$ with a compactly supported smooth function and use Young's convolution inequality and Theorem \ref{thm_main_Lp_moment} with $p=p_\alpha$.

For $\e>0$, let 
   \begin{align*}
       E_\e:=\bigcup_{i=1}^{M_\e} \Gamma_{C\e, t_i}.
    \end{align*}
We can choose a sufficiently large $C$ such that the support of $u_\e$ is contained in $E_\e$. Note that $|E_\e| \lesssim \e^{\frac{d^2+d-2\alpha}{2}}$. Let $\psi$ be a compactly supported function in $\R^d$. Applying Lemmas \ref{lemma_1} and \ref{lemma_2}, we have
\begin{align*}
    |\langle u, \psi\rangle|^2
    &=\lim\limits_{\e\to 0}|\langle u_\e, \psi\rangle|^2\\
    &\leq \lim\limits_{\e\to 0} \Vert u_\e\Vert_2^2\int_{E_{\e}}|\psi|^2\\
    &\lesssim\Vert \psi\Vert_\infty^2\lim_{\e\to 0}\e^{-\frac{d^2+d-2\alpha}{2}}|E_{\e}|\sum_{j=-\infty}^\infty a_jb_{j,\e}\\
    &=\Vert \psi\Vert_\infty^2\lim_{\e\to 0}\sum_{j=-\infty}^\infty a_jb_{j,\e}\\
    &=0.
\end{align*}
In the last equality, we used the dominated convergence theorem. Since this holds for any compactly supported smooth function $\psi$, we conclude that $u \equiv 0$. Thus, $f\equiv 0$ as desired.

\end{proof}

\subsection{ \texorpdfstring{$L^2$}{L2} estimate for  \texorpdfstring{$u_\e$}{u epsilon}}
In this section, we will prove Lemma \ref{lemma_1}. First, we will show that
\begin{equation}\label{lem1_eq1}
    \norm{u_\e}_2^2 \lesssim_N \e^{-\frac{d^2+d-2\alpha}{2}} \sum_{j \in \Z} a_j b_{j,\e}
\end{equation}
for $ 0 < \e <1$, where
$$
a_j=\min\{2^{-jN},1\} 2^{jd\cdot \frac{p_\alpha-2}{p_\alpha}} \qquad \text{and} \qquad b_{j,\e}=\bigg(\int_{A_{\e,j}}\bigg(\sum\limits_i|f_{\e,i}(x)|^{2}\bigg)^{\frac{p_\alpha}{2}} dx\bigg)^{\frac{2}{p_\alpha}}.    
$$
The region $A_{\e,r}$ will be specified later. Then, we will show that $a_j$ and $b_{j,\e}$ satisfy the desired properties.

\begin{proof}[Proof of \eqref{lem1_eq1}]
    Observe that by the construction, for each $i$, the support of $\widehat{f_{\e,i}}\ast \eta_{\e,i}$ is contained in $ CB_{\e,i}$ for some large constant $C$ which does not depend on $\e$ and $i$. Hence, the supports of  $\widehat{f_{\e,i}}\ast \eta_{\e,i}$ have bounded overlaps. By Plancherel theorem, we have
    \begin{align}\label{sum1}
        \norm{u_\e}_2^2
        &=\int \left|\sum\limits_i \widehat{f_{\e,i}}\ast \eta_{\e,i}(\xi)\right|^2d\xi \nonumber\\
        &\lesssim\sum\limits_i\int \big| \widehat{f_{\e,i}}\ast \eta_{\e,i}(\xi)\big|^2d\xi\nonumber \\
        &=\sum\limits_i\int |f_{\e,i}(x)|^2|\widehat{\eta_{\e,i}}(x)\big|^2dx.
    \end{align}
    Recall that $\widetilde{B}_{\e,i}$ is the dual rectangular box of $B_{\e,i}$ (See Definition \ref{dualbox}). For each $1\leq i\leq M_\e$ and each $j\in \mathbb{Z}$, we let $A_{\e,i,j}=2^j\widetilde{B}_{\e,i} \setminus 2^{j-1}\widetilde{B}_{\e,i}$. Note that $\widehat{\eta_{\e,i}}$ decays rapidly outside of the dual box $\widetilde{B}_{\e,i}$ with side lengths $(\e^{-1},\dots, \e^{-d})$. Hence, 
    \eqref{sum1} is bounded by
    \begin{align*}
        &\lesssim\sum\limits_i\sum_{j\in \Z} \int_{A_{\e,i,j}} |f_{\e,i}(x)|^2|\widehat{\eta_{\e,i}}(x)\big|^2dx\\
        &\lesssim\sum\limits_i\sum_{j\in \Z} \sup\limits_{y\in A_{\e,i,j}} \big|\widehat{\eta_{\e,i}}(y)\big|^2\int_{A_{\e,i,j}} |f_{\e,i}(x)|^2dx\\
        &\lesssim_N\sum\limits_i\sum_{j\in \Z} \min\{2^{-jN},1\}\int_{A_{\e,i,j}} |f_{\e,i}(x)|^2dx,
    \end{align*}
    with $N\in \N$ large enough. 
    
    For each $j\in\Z$, we let
    $$
        A_{\e,j}=\bigcup_{i=1}^{M_\e} A_{\e,i,j}.
    $$
    Invoking H\"older inequality with the exponent $p_\alpha/2$, $\Vert u_\e\Vert_2^2$ is dominated by
    \begin{align}
        &\sum\limits_{i}\sum_{j\in \Z} \min\{2^{-jN},1\}\int_{A_{\e,j}} |f_{\e,i}(x)|^2dx \nonumber \\
        &=\sum_{j\in \Z}\min\{2^{-jN},1\} \int_{{A_{\e,j}}}\sum\limits_i|f_{\e,i}(x)|^2dx \nonumber \\
        &\lesssim \sum_{j\in \Z}\min\{2^{-jN},1\} \bigg(\int_{A_{\e,j}}\bigg(\sum\limits_i|f_{\e,i}(x)|^{2}\bigg)^{\frac{p_\alpha}{2}} dx\bigg)^{\frac{2}{p_\alpha}} |{A_{\e,j}}|^{1-\frac{2}{p_\alpha}}. \label{eq_u_epsilon}
        \end{align}
    Since $M_\e\lesssim\e^{-\alpha}$, we have 
    \begin{align*}
        |{A_{\e,j}}|\lesssim \e^{-\alpha}\prod_{l=1}^d2^j(\e^{-1})^{l}=\e^{-\frac{d^2+d}{2}-\alpha}2^{jd},
    \end{align*}
    which implies that
    \begin{align*}
        |{A_{\e,j}}|^{1-\frac{2}{p_\alpha}}\lesssim \e^{-\frac{d^2+d-2\alpha}{2}}2^{jd\cdot \frac{p_\alpha-2}{p_\alpha}}.
    \end{align*}
    Combining with \eqref{eq_u_epsilon}, we obtain that
    \begin{align*}
        \Vert u_\e\Vert_2^2 
        &\lesssim_N \e^{-\frac{d^2+d-2\alpha}{2}} \sum_{j\in \Z}\min\{2^{-jN},1\} 2^{jd\cdot \frac{p_\alpha-2}{p_\alpha}}\bigg(\int_{A_{\e,j}}\bigg(\sum\limits_i|f_{\e,i}(x)|^{2}\bigg)^{\frac{p_\alpha}{2}} dx\bigg)^{\frac{2}{p_\alpha}}\\
        &=\e^{-\frac{d^2+d-2\alpha}{2}}\sum\limits_{j\in \Z}a_j b_{j,\e},
    \end{align*}
\end{proof}
In order to prove the properties of $b_{j,\e}$, we need to show the following.

    \begin{lemma}\label{lemma_2.5}
        If $p \geq 2$, we have
        \begin{equation*}
            \norm{\bigg(\sum_{i=1}^{M_\e} |f_{\e,i}|^2\bigg)^{1/2}}_{p} \lesssim \norm{f}_{p}
        \end{equation*}
        uniformly in $0 < \e< 1$.
    \end{lemma}
    
    \begin{lemma}\label{lemma_3}
        Assume that $f \in L^p(\R^d)$ for some $p >2$. For fixed $j \in \Z$, we have
        \begin{equation*}
            \bigg(\int_{|x| \geq 2^j \e^{-1} }\bigg(\sum\limits_i|f_{\e,i}(x)|^{2}\bigg)^{\frac{p}{2}} dx\bigg)^{\frac{1}{p}} \longrightarrow 0, \qquad \text{as}\  \e \rightarrow 0.
        \end{equation*}
    \end{lemma}

    Lemma \ref{lemma_2.5} is a result of Rubio de Francia's inequality, which is Theorem 1.2 in \cite{Rubio85}.
    \begin{proof}[Proof of Lemma \ref{lemma_2.5}]
        Observe that the Fourier transform of each $f_{\e,i}$ in the first variable is supported in an interval of length $\sim \e$, for all $1\leq i\leq M_\e$. In \cite{Rubio85}, the proof was reduced to the case of the smooth operator $G$ associated with a sequence of finite overlapping intervals (See Definition 3.3 in \cite{Rubio85}). Thus, we can apply Rubio de Francia's result to our setting in the first variable. Then, we get
\begin{align*}
   \int \bigg(\sum_i|f_{\e,i}(x_1,x')|^2\bigg)^\frac{p}{2} dx_1 \leq C_p \int |f(x_1,x')|^p dx_1, \qquad \forall x'\in \R^{d-1}.
\end{align*}

Integrating in the variables $x'$, we obtain that
\begin{align*}
\bigg\Vert\bigg(\sum_{i=1}^{M_\e} |f_{\e,i}|^2\bigg)^{1/2}\bigg\Vert_{p}\lesssim ||f||_{p}.
\end{align*}
    \end{proof}

    To prove Lemma \ref{lemma_3}, we will use the weighted version of Rubio de Francia's inequality, which is Theorem 6.1 in \cite{Rubio85}. Recall that a \textit{weight} is a nonnegative locally integrable function on $\R^d$ that takes values in $(0,\infty)$ almost everywhere.
\begin{definition}
    Let $1<p<\infty$. A weight $w$ is said to be of class $A_p$ if the $A_p$ Muckenhoupt characteristic constant $[w]_{A_p}$ of $w$ is finite, i.e.
    \begin{align*}
        [w]_{A_p}=\sup_{Q \text{ cubes in }\R^d}\bigg(\frac{1}{|Q|}\int_Qw(x)dx\bigg)\bigg(\frac{1}{|Q|}\int_Qw(x)^{-\frac{1}{p-1}}dx\bigg)^{p-1}<\infty.
    \end{align*}
\end{definition}
For the properties of the Muckenhoupt weight, see Chapter 7 in \cite{Grafakosbook}.

    \begin{proof}[Proof of Lemma \ref{lemma_3}] 
    Fix $j\in \Z$, and let $\delta$ be a number such that $0< \delta<\frac{p}{2}-1$. For $x\in \R^d$, we write $x=(x_1,x')$ where $x'=(x_2,\dots, x_d)\in \R^{d-1}$. Let $w:\R^d\to (0,\infty)$ be a weight in the variable $x_1$ defined by
    \begin{align*}
        w(x_1,x') = \left\{
    \begin{matrix}
        \left( \frac{|x_1|}{ 2^j \e^{-1}} \right)^{p/2 - 1 -\delta } & \text{if} \ |x_1 | \leq 2^j \e^{-1} \ \text{and} \ |x'| \leq 2^{j} \e^{-1},\\
        1 & \text{otherwise}.
    \end{matrix}\right.
    \end{align*}
    We can show that $[w(\cdot,x')]_{A_{p/2}}\lesssim 1$ uniformly in $x'$, $j$, and $\e$. Indeed, if $|x'|\geq 2^j\e^{-1}$, then $w(x_1,x')=1$ and $[w(\cdot,x_1)]_{A_{p/2}}=1$. If $|x'|\leq 2^j\e^{-1}$, then
    \begin{align*}
        w(x_1,x')=\min \bigg\{ \left(\frac{|x_1|}{ 2^j \e^{-1}} \right)^{p/2 - 1 -\delta },1\bigg\}.
    \end{align*}
    Observe that $w_1(x_1)=|x_1|^{p/2 - 1 -\delta }$ is an $A_{p/2}$ weight, and $[w_1(\lambda x_1)]_{A_{p/2}} = [w_1]_{A_{p/2}}$ for any $\lambda>0$. For example, see Proposition 7.1.5 and Example 7.1.7 in \cite{Grafakosbook}. Since $\delta < p/2-1$, $w(\cdot, x')$ is also an $A_{p/2}$ weight, with
    \begin{align*}
        [w(\cdot,x')]_{A_{p/2}}=[\min\{ w_1,1\}]_{A_{p/2}}\lesssim_p [w_1]_{A_{p/2}} \lesssim 1.
    \end{align*}

    Recall that the Fourier transform of $f_{\e,i}$ in the first variable is supported in an interval of length $\sim \e$, for all $1\leq i\leq M_\e$, and those intervals overlap at most finitely many times. Since $p > 2$ and $[w (\cdot, x')]_{A_{p/2}} \lesssim 1$ uniformly in $x'$, applying the weighted version of Rubio de Francia's inequality, we obtain that 
    \[
        \int \bigg(\sum\limits_{i=1}^{M_\e}|f_{\e,i}(x_1, x')|^{2}\bigg)^{\frac{p}{2}} w(x_1, x') dx_1\lesssim  \int |\sum_if_{\e,i}(x_1,x')|^{p} w(x_1, x') dx_1.
        \] 
    Integrating the inequality above with respect to $x'$, we get
    \[
        \bigg(\int_{|x| \gtrsim 2^j \e^{-1} }\bigg(\sum\limits_i|f_{\e,i}(x)|^{2}\bigg)^{\frac{p}{2}} dx\bigg)^{\frac{1}{p}}\lesssim \norm{\left( \sum_i |f_{\e,i} |^2\right)^{1/2}}_{L^{p}(w)} \lesssim \norm{f}_{L^{p}(w)}.
    \]
    Now, we have 
    \begin{align*}
        \norm{f}_{L^{p}(w)}^{p}
    &= \int_{|x| \geq (2^j \e^{-1})^{1/2}} |f|^{p}w(x)dx+\int_{|x| \leq (2^j \e^{-1})^{1/2}} |f|^{p}w(x)dx\\
    &\leq \int_{|x| \geq (2^j \e^{-1})^{1/2}} |f|^{p}dx + (2^{-j}\e)^{\frac{1}{2}(\frac{p}{2}-1-\delta)} \int |f|^{p} dx.
    \end{align*}
    Since $\norm{f}_p < \infty$, we establish the result by letting $\e \rightarrow 0$.      
    \end{proof}
    Now, we are ready to prove Lemma \ref{lemma_1}.
    \begin{proof}[Proof of Lemma \ref{lemma_1}]
    We have established \eqref{lem1_eq1}. Thus, it suffices to show that $a_j$ and $b_j$ satisfies the desired properties. For sufficiently large $N$, we have $\sum_{j} |a_j| < \infty$. By Lemma \ref{lemma_2.5}, we get $|b_{j,\e} | \lesssim \norm{f}_{p_\alpha}^2$. Now we observe that the ball of radius $\sim 2^j \e^{-1}$ centered at the origin is contained in the intersection of dual boxes $\widetilde{B}_{\e,i}$. Thus for each $j\in \Z$, we have
 \begin{align*}
     {A}_{\e,j}\subseteq \R^d\setminus B(0,2^{j-2}\e^{-1}).
 \end{align*}
 Therefore, Lemma \ref{lemma_3} implies that $b_{j,\e} \rightarrow 0$ as $\e \rightarrow 0$.
    \end{proof}

\subsection{Weak convergence of \texorpdfstring{$\{u_\e\}_{\e>0}$}{u epsilon}}
In this section we will prove Lemma \ref{lemma_2}.
\begin{proof}[Proof of Lemma \ref{lemma_2}]
We will show that 
\[ \lim_{\e \rightarrow 0} \langle u - u_\e , \psi \rangle =0. \]
Note that by Plancherel's theorem, we can write 
\begin{align*}
    |\langle u - u_\e , \psi \rangle|
    &=|\langle \widehat{f}-\sum_i\widehat{f_{\e,i}}\ast\eta_{\e,i} , \psi \rangle|=\big|\langle \sum_if_{\e,i}(1
    -\widehat{\eta_{\e,i}} ), \widehat{\psi} \rangle\big|.
\end{align*}

Thus, we have
\begin{align*}
    |\langle u - u_\e , \psi \rangle|
    &\leq \int_{\{|x|\leq \e^{-\frac{1}{4d}}\}} \bigg(\sum_i |f_{\e,i}(x) (\widehat{\eta_{\e,i}}(x) -1)| \bigg)|\widehat{\psi}(x)|dx\\
    &\hspace{1cm}+\int_{\{|x|>\e^{-\frac{1}{4d}}\}}\bigg(\sum_i |f_{\e,i}(x) (\widehat{\eta_{\e,i}}(x) -1)| \bigg)|\widehat{\psi}(x)|dx\\
    &=I+II.
\end{align*}
First, we estimate $I$. By Cauchy-Schwarz inequality, we can bound $I$ by
\begin{align}\label{eq_estimate_I}
    I\lesssim \int_{\{|x|\leq \e^{-\frac{1}{4d}}\}}\bigg(\sum_i|f_{\e,i}(x)|^2\bigg)^{1/2}\bigg(\sum_i|\widehat{\eta_{\e,i}}(x)-1|^2\bigg)^{1/2}|\widehat{\psi}(x)|dx.
\end{align}
Observe that $\widehat{\eta_{\e,i}}(0)=1$. Since $|x|\leq \e^{-\frac{1}{4d}}$, we have that for each $1\leq i\leq M_\e$,
\begin{align*}
    |\widehat{\eta_{\e,i}}(x)-1| =|\widehat{\eta_{\e,i}}(x)-\widehat{\eta_{\e,i}}(0)| \lesssim \e^{1- \frac{1}{4d}}.
\end{align*}
Therefore, using the fact that $M_\e\lesssim \e^{-\alpha } \leq  \e^{-1}$, we have
\begin{align*}
\bigg(\sum_{i=1}^{M_\e}|\widehat{\eta_{\e,i}}(x)-1|^2\bigg)^{1/2}\lesssim \e^{\frac{1}{2} -\frac{1}{4d}}.
\end{align*}

We plug this into \eqref{eq_estimate_I}, and then apply Holder inequality with exponent $p$. Then we get
\begin{align*}
    I&\lesssim \e^{\frac{1}{2} -\frac{1}{4d}}||\widehat{\psi}||_{\infty}\,\int_{\{|x|\leq \e^{-\frac{1}{4d}}\}}\bigg(\sum_i|f_{\e,i}(x)|^2\bigg)^{1/2}dx\\
    &\lesssim \e^{\frac{1}{2} -\frac{1}{4d}}||\widehat{\psi}||_{\infty}\,\bigg(\int \bigg(\sum_i|f_{\e,i}(x)|^2\bigg)^{p/2}dx\bigg)^{1/p}\cdot \bigg(\int_{\{|x|\leq \e^{-\frac{1}{4d}}\}} 1 dx\bigg)^{1-\frac{1}{p}}\\
    &=\e^{\frac{1}{4} ( 1-\frac{1}{d}+\frac{1}{p})}||\widehat{\psi}||_{\infty}\,\bigg\Vert\bigg(\sum_{i=1}^{M_\e} |f_{\e,i}|^2\bigg)^{1/2}\bigg\Vert_{p}.
\end{align*}
By Lemma \ref{lemma_2.5}, we obtain that
\begin{align}\label{eq_estimate_I_result}
   I \lesssim  \e^{\frac{1}{4} ( 1-\frac{1}{d}+\frac{1}{p})} \, ||\widehat{\psi}||_{\infty}\,||f||_{p}.
\end{align}

Now we bound the second integral. Applying Holder inequality with the exponent $p$, we have
\begin{align*}
    II
    &\lesssim \sum_{i=1}^{M_\e} \int_{\{|x|>\e^{-\frac{1}{4d}}\}} |f_{\e,i}(x)||\widehat{\psi}(x)|dx\\
    &\lesssim \sum_{i=1}^{M_\e} \bigg(\int |f_{\e,i}(x)|^{p}dx\bigg)^{1/p}\cdot \bigg(\int_{\{|x|>\e^{-\frac{1}{4d}}\}} |\widehat{\psi}(x)|^{p'}dx\bigg)^{1/p'}.
\end{align*}
Since $\widehat{\psi}$ is a Schwartz function, we have
\begin{align*}
    |\widehat{\psi}(x)|\lesssim (1+|x|)^{-8d},\qquad \forall x\in \R^d.
\end{align*}
Also, Lemma \ref{lemma_2.5} implies that  $||f_{\e,i}||_{L^{p}(\R^d)}\lesssim||f||_{L^{p}(\R^d)}$ uniformly in $\e$ and $i$. Thus, we can bound $II$ by
\begin{align}\label{eq_estimate_II_result}
        II&\lesssim \sum_{i=1}^{M_\e} ||f_{\e,i}||_{p}\cdot \bigg(\int_{\{|x|>\e^{-\frac{1}{4d}}\}} |x|^{-8d p'}dx\bigg)^{1/p'} \nonumber \\
    &\lesssim \sum_{i=1}^{M_\e} ||f_{\e,i}||_{p} \e^{2- \frac{d-1}{4d}(1-\frac{1}{p})} \nonumber \\
    &\lesssim  \e^{1- \frac{d-1}{4d}(1-\frac{1}{p})} \,||f||_{p} \nonumber \\
    &\lesssim  \e^{\frac{1}{4}(3+\frac{1}{p})} \,||f||_{p}
\end{align}
Combining \eqref{eq_estimate_I_result} and \eqref{eq_estimate_II_result}, we have
\begin{align*}
    \langle u-u_\e,\psi\rangle \lesssim  \e^{\frac{1}{4} ( 1-\frac{1}{d}+\frac{1}{p})}\,\max(||\widehat{\psi}||_{\infty},1)\,||f||_{p}.
\end{align*}
Thus, $\lim\limits_{\e\to 0} \langle u-u_\e,\psi\rangle=0$ as desired.
\end{proof}

\section{Optimality of the range of \texorpdfstring{$p$}{p} in \texorpdfstring{$L^p$}{Lp}-integrability}\label{sec3}
We will prove Theorem \ref{thm_main_Lp_optimal} by a probabilistic method. An example with $\alpha =1$ was considered in Theorem 1.3 in \cite{ACK_book}. Thus, we assume that $0 < \alpha <1$ throughout this section. We abbreviate almost surely to a.s. and we will simply write $\gamma(t)$ to denote $\gamma_d(t)$. 

Let us consider a nondecreasing sequence $\{m_j\}_{j \in \N}$ such that $m_j \geq 2$ and
\begin{equation*}
    m_{j+1} \lesssim_\e (m_1 m_2 \cdots m_j)^{\e}, \qquad \forall \e >0
\end{equation*}
For each $j$, we write $M_j = m_1 m_2 \cdots m_j$. The collection of $M_j^{-1}$-intervals $\cI_j$ is defined by
\[
\{ M_j^{-1} (m + [0,1)) : 0 \leq m \leq M_j-1   \}.
\]
We consider a sequence of random functions $\mu_j$ which satisfies the following conditions for some deterministic nondecreasing sequence $\{\beta_j \}_{j \in \N}$:
\begin{itemize}
    \item $\mu_0 = \mathbf{1}_{[0,1]}$.
    \item $\mu_j = \beta_j \mathbf{1}_{E_j}$ where $E_j$ is a union of intervals in $\cI_j$.
    \item $\mathbb{E}(\mu_{j+1} (x) |E_j) = \mu_j(x)$ for all $x \in [0,1]$.
    \item The sets $I_j \cap E_{j+1}$ are chosen independently for each $I_j \in \mathcal{I}_j$ conditioned on $E_j$.
\end{itemize}
We identify the functions $\mu_j$ with the measures $\mu_j dx$ and denote that $\norm{\mu_j}  = \mu_j([0,1])$.
Let $\alpha$ be a number such that 
\begin{equation}\label{cond_b_jM_j}
    1-\alpha = \lim_{j \rightarrow \infty} \frac{\log(\beta_j)}{\log(M_j)}.
\end{equation}
Note that it implies that $M_j^{1-\alpha-\e} \lesssim_\e \beta_j \lesssim_\e M_j^{1-\alpha +\e}$.

Recall that $\gamma (t) :=\gamma_d(t) = (t, t^2, \cdots, t^d)$ and $t \in [0,1]$. We define a measure $\nu_j$ on $\Gamma_d$ by
\[
\int f(x_1, x_2, \cdots, x_d) d\nu_j =  \int f(\gamma(t)) d\mu_j(t).
\]
Then a.s. the measure $\mu_j$ converges weakly to a measure $\mu$ supported on a subset of $[0,1]$ and a.s. $\nu_j$ converges weakly to a measure $\nu $ supported on a subset of $\Gamma_d$. For more properties of the measure $\mu$, see \cite{ShSu17}. 

\begin{prop}\label{Prop_optimal}
    If $p > p_\alpha$, then a.s. $\norm{\widehat{\nu}}_p \lesssim 1$.
\end{prop}
Throughout this section, when we say an inequality holds a.s., the corresponding implicit constant may depend on the limiting measure $\nu$. However, the probability that the implicit constant exists is $1$.

Theorem \ref{thm_main_Lp_optimal} follows from Proposition \ref{Prop_optimal}. 
\begin{proof}[Proof of Theorem \ref{thm_main_Lp_optimal} using Proposition \ref{Prop_optimal}]
For sufficiently large number $m \in \N $, let $m_j =m$ for each $j \in \N$ so that $M_j = m^j$ and let $\beta_j = M_j^{1-\alpha}$.

For each $j$, we choose a subset $S_j$ of $\{ 1, 2, \cdots, m\}$ of size $\sim m^\alpha$ such that $|S_{j_1}| \cdot |S_{j_1+1}|\cdots |S_{j_2}| \sim m^{(j_2-j_1+1)\alpha}$ for any $j_2 >j_1$. We choose each $I_j \cap E_{j+1}$ by random translation of $M_{j+1}^{-1}(S_{j+1}+[0,1]))$ as in Section 8 of \cite{LW18}. More specifically, we randomly choose a number $n(I_j)$ such that $\mathbb{P}(n(I_j) = n) =m^{-1}$ for each $I_j$ and $1 \leq n \leq m$. We let 
    $$S_{j+1,n(I_j)} = S_{j+1}+n(I_j) \mod m.$$
    Then, $I_j \cap E_{j+1} = \ell(I_j) + M_{j+1}^{-1} (S_{j+1,n(I_j)} + [0,1))$ where $\ell(I_j)$ is the left endpoint of $I_j$. Note that $\nu$ is an Ahlfors-David regular $\alpha$-measure, i.e. $ \nu(B(x,r))) \sim r^\alpha $ for all $ x \in \spt \nu$ and $ 0< r < \diam(\spt \nu)$. Thus, $N(\spt {\nu}, \e) \sim \e^{-\alpha}$ for all $ 0< \e <1$.

    By Proposition \ref{Prop_optimal}, there exists $\nu$ such that $\norm{\widehat{\nu}}_p \lesssim 1$. Therefore, $\widehat{\nu}$ is the desired function $f$.
\end{proof}

When $d=2$, an estimate on the Fourier decay of the measure $\nu$ was studied in \cite{Ryou23}.  

\begin{proof}[Proof of Proposition \ref{Prop_optimal} when $d =2$]

    By (2) in Theorem 1.2 \cite{Ryou23}, a.s. we have 
    \[
    |\widehat{\nu}(\xi)| \lesssim_\e (1+|\xi|)^{-\alpha/2+\e}, \qquad \forall \xi \in \mathbb{R}^2.
    \]
    Hence, $\norm{\widehat{\nu}}_p \lesssim 1$ for $p > p_\alpha$.
\end{proof}
Let us turn to the proof of Proposition \ref{Prop_optimal} in higher dimensions.

First, we will consider when $d \geq 4$. For $s_1\geq 0$ and $1 \leq 2^{s_2} \leq M_j2^{s_1}$, we define the sets 
\[
\Omega_{j} (s_1):=\{\xi \in \R^d : \exists t \in E_{j} \ \text{such that} \ |\langle \gamma^{(k)} (t), \xi \rangle | \leq (M_j2^{s_1})^k, \ \forall 1 \leq k \leq d \},
\]
and
\begin{align*}
    \Omega_{j} (s_1,s_2):=\{\xi \in \Omega_j(s_1) : \exists t \in E_{j} \ \text{such that} \ |\langle \gamma^{(k)} (t), \xi \rangle | \leq (M_j2^{s_1}) 2^{s_2(k-1)}, \ \forall 2 \leq k \leq d \}.
\end{align*}

\begin{lemma}\label{lem_mup_meas}
    For any $\e >0$, a.s., we have the following estimates.
    \begin{itemize}
        \item[(i)] If $1\leq 2^{s_2 } \leq (2^{s_1} M_j)^{\frac{1}{d-1}}$, then
    \begin{equation}\label{eq_Omega_meas1}
        |\Omega_j(s_1,s_2) | \lesssim_\e M_j^\e (M_j 2^{s_1})^d 2^{s_2(\frac{d(d-1)}{2}+\alpha d)} 2^{s_1\frac{d(1-\alpha)}{d-1}}.
        \end{equation}
        \item[(ii)] If $(2^{s_1} M_j)^{\frac{1}{d-1}} \leq 2^{s_2 } \leq 2^{s_1 } M_j $, then
    \begin{equation}\label{eq_Omega_meas2}
        |\Omega_j(s_1,s_2)| \lesssim_\e M_j^\e (M_j 2^{s_1})^d 2^{s_2 \frac{d(d-1)}{2} } M_j^\alpha 2^{s_1} ( M_j 2^{s_1-s_2})^{\frac{\alpha}{d-2}}2^{s_1\frac{1-\alpha}{d-1}}.
    \end{equation}
    \end{itemize}
\end{lemma}
To estimate $|\widehat{\nu}|$, we define the set $\widetilde{\Omega}_{j+1}(s_1,s_2)$ as follows: 
\begin{align*}
	\widetilde{\Omega}_{j+1}(s_1,s_2) &=
	\Omega_{j+1}(s_1,s_2) \backslash (\Omega_{j+1}(s_1,s_2-1) \cup \Omega_{j+1}(s_1-1)) \quad \text{if } s_1, s_2 \neq 0,\\
	\widetilde{\Omega}_{j+1}(s_1,0) &=
	\Omega_{j+1}(s_1,0) \backslash \Omega_{j+1}(s_1-1) \quad \text{if } s_1 \neq  0,\\
	\widetilde{\Omega}_{j+1}(0,s_2) &=
	\Omega_{j+1}(0,s_2) \backslash \Omega_{j+1}(0,s_2-1) \quad \text{if } s_2 \neq 0,\\
	\widetilde{\Omega}_{j+1}(0,0) &=
	\Omega_{j+1}(0,0).\\
\end{align*}
\begin{lemma}\label{lem_mup_osci}
    If $\xi \in \widetilde{\Omega}_{j+1}(s_1,s_2)$, for any $\e >0$, a.s. we have 
    \begin{equation}\label{Omega_osci}
        |\widehat{\nu}_{j+1}(\xi) - \widehat{\nu}_j(\xi) | \lesssim_\e (M_j 2^{s_1})^\e
 M_{j}^{-\frac{\alpha}{2}} 2^{-s_1(1-\frac{\alpha }{2})-s_2\frac{\alpha}{2}}.
    \end{equation}
\end{lemma}
\begin{remark}
    When $\alpha=1$, for each fixed $\xi \in \R^d$, it suffices to consider one oscillatory integral
    $$
    \int_0^1 e(-\langle \xi, \gamma(t)\rangle)dt.
    $$
    However, when $0 < \alpha <1$, $\widehat{\nu}_j (\xi)$ is a summation of oscillatory integrals of form
    $$
    \beta_{j}\int_{I_j}  e(-\langle \xi, \gamma(t)\rangle)dt
    $$
    where $I_j \in \mathcal{I}_j $ and $I_j \subseteq E_j $. Therefore, we need to show that the summation of oscillatory integrals has a good upper bound and this is why we decompose $\Omega_j(s_1)$ further into $\widetilde{\Omega}_j(s_1,s_2)$.
\end{remark}

We prove Proposition \ref{Prop_optimal} using Lemmas \ref{lem_mup_meas} and \ref{lem_mup_osci} first. The proof of lemmas will be given later.

\begin{proof}[Proof of Proposition \ref{Prop_optimal} when $d \geq 4$]
    It suffices to show that $\norm{\widehat{\nu}_n}_p \lesssim 1 $ uniformly in  $n \in \N$. By letting $n \rightarrow \infty$, we have $\norm{\widehat{\nu}}_p \lesssim 1$. By triangle inequality, we get
    \[
        \norm{\widehat{\nu}_n }_p \lesssim \sum_{ 1 \leq j \leq n-1 } \norm{\widehat{\nu}_{j+1} - \widehat{\nu}_j }_p  +1.
    \]
    For each $j$, we have
    \begin{align}\label{eq_v_j}
        \norm{\widehat{\nu}_{j+1} - \widehat{\nu}_j }_p^p 
        &\leq 
        \sum_{s_1\geq 0}\bigg(\sum_{ 1 \leq 2^{s_2} \leq (2^{s_1}M_{j+1})^{\frac{1}{d-1}} }\int_{\widetilde{\Omega}_{j+1}(s_1,s_2)} |\widehat{\nu}_{j+1} - \widehat{\nu}_j|^p d\xi\bigg) \nonumber\\
        &\hspace{1cm}+ 
         \sum_{s_1\geq 0}\bigg(\sum_{ (2^{s_1}M_{j+1})^{\frac{1}{d-1}} \leq 2^{s_2} \leq 2^{s_1}M_{j+1}}\int_{\widetilde{\Omega}_{j+1}(s_1,s_2)} |\widehat{\nu}_{j+1} - \widehat{\nu}_j|^p d\xi\bigg) \nonumber\\
        &= \mathcal{J}_1+\mathcal{J}_2.
    \end{align}
    First, we estimate the sum $\mathcal{J}_1$. Recall that $m_{j+1} \lesssim_\e M_j^\e$. For fixed $s_1\geq 0$, since $1 \leq 2^{s_2} \leq (2^{s_1} M_{j+1})^{\frac{1}{d-1}}$,  Lemmas \ref{lem_mup_meas} and \ref{lem_mup_osci} imply that each integral in the sum $\mathcal{J}_1$ can be bounded by
    \begin{align*}
        (M_j 2^{s_1})^\e M_j^{d-\frac{\alpha p}{2}}2^{s_1(\frac{\alpha-2}{2}p+d+\frac{d(1-\alpha)}{d-1})}2^{s_2(-\frac{\alpha p}{2}+\frac{d(d-1)}{2}+\alpha d)}.
    \end{align*}
    Let $p_2(\mathcal{J}_1)$ denote the exponent of $2^{s_2}$, i.e., $p_2(\mathcal{J}_1) = -\frac{\alpha p}{2}+\frac{d(d-1)}{2}+\alpha d$. Then, we can check that the sum over $1 \leq 2^{s_2} \leq (2^{s_1} M_{j+1})^{\frac{1}{d-1}}$ is bounded by
    \begin{align}\label{sum_s2_1}
            &\lesssim_\e \begin{cases}
             (M_j 2^{s_1})^\e M_{j}^{ \frac{d}{2(d-1)}(-\alpha p +3(d-1) +2\alpha)}2^{s_1 p_1( \mathcal{J}_1)}&\text{ if } p_2(\mathcal{J}_1)> 0,\\
             (M_j 2^{s_1})^\e M_{j}^{d-\frac{\alpha p}{2}}2^{s_1 p_1( \mathcal{J}_1)}  &\text{ if } p_2(\mathcal{J}_1)\leq 0,
        \end{cases}
    \end{align}
    where 
    \begin{align*}
            p_1(\mathcal{J}_1) &= \begin{cases}
            \frac{\alpha-2}{2}p + \frac{3d}{2} + \frac{1}{(d-1)}(-\frac{\alpha p}{2} +d)&\text{ if } p_2(\mathcal{J}_1)> 0,\\
            \frac{\alpha-2}{2}p+d+\frac{d(1-\alpha)}{d-1}  &\text{ if } p_2(\mathcal{J}_1) \leq 0.
        \end{cases}
    \end{align*}
    In either case, we will show that $p_1(\mathcal{J}_1) <0$.  When $p_2(\mathcal{J}_1) > 0$, we can rewrite  $p_1( \mathcal{J}_1)$ as
    \[
    \begin{split}
    	\frac{d}{2(d-1)}(-\alpha p + 3(d-1)+2\alpha) + \frac{1-\alpha}{d-1} (-p(d-1)+d)={\Romannum{1}}+ {\Romannum{2}}.
    \end{split}
    \]
    It is easy to see that ${\Romannum{2}}$ is negative. For ${\Romannum{1}}$, by the assumptions that $d \geq 4$, $p >p_\alpha$, and $0<\alpha<1$, we can check that
    \begin{equation}\label{eq_rom2}
    	{\Romannum{1}}< \frac{d}{4(d-1)} (-d^2+5d-6+2\alpha)\leq -\frac{1}{4}d(d-4) \leq 0.    
    \end{equation}
    When $p_2(\mathcal{J}_1) \leq 0$, using the assumptions that $\alpha <1$ and $p > p_\alpha $, we obtain that
    \[
    	p_1(\mathcal{J}) = -p+d+ \frac{d}{d-1} + \alpha \left(\frac{p}{2} - \frac{d}{d-1}\right) < -\frac{p}{2} +d  < 0.
    \]
    Since $p_1(\mathcal{J}_1)<0$, we have
    \begin{align}\label{eq_sum_J1}
        \mathcal{J}_1\lesssim_\e M_{j}^{\e-\frac{\alpha p}{2}+d}+M_{j}^{\e+ \frac{d}{2(d-1)}(-\alpha p +3(d-1) +2\alpha)}.
    \end{align}
    To estimate the sum $\mathcal{J}_2$, we fix $s_1\geq 0$. Since $ (2^{s_1} M_{j+1})^{\frac{1}{d-1}} \leq 2^{s_2} \leq 2^{s_1} M_{j+1}$, we apply Lemmas \ref{lem_mup_meas} and \ref{lem_mup_osci} to bound each integral in the sum by a constant multiple of
    \begin{align}\label{eq_est_sumJ2}
         (M_j 2^{s_1})^\e M_{j}^{-\frac{\alpha p}{2}+d+\frac{\alpha (d-1)}{(d-2)}}2^{s_1(\frac{(\alpha-2)p}{2}+d+1+\frac{\alpha}{d-2}+\frac{1-\alpha}{d-1})}2^{s_2 p_s(\mathcal{J}_2)}.
    \end{align}
    where $p_2(\mathcal{J}_2) = -\frac{\alpha p}{2}+\frac{d(d-1)}{2}-\frac{\alpha}{d-2}$. The sum over $s_2$ with $ (2^{s_1} M_{j+1})^{\frac{1}{d-1}} \leq 2^{s_2} \leq 2^{s_1} M_{j+1}$ is bounded by
    \begin{align}\label{sum_s2_2}
        &\lesssim_\e \begin{cases}
             (M_j 2^{s_1})^\e M_{j}^{-\alpha p+\frac{d(d+1)}{2}+\alpha}2^{s_1 p_1(\mathcal{J}_2)}  &\text{ if } p_2(\mathcal{J}_2)>0,\\
            (M_j 2^{s_1})^\e M_{j}^{ \frac{d}{2(d-1)}(-\alpha p +3(d-1) +2\alpha)}2^{s_1 p_1(\mathcal{J}_2)}
           &\text{ if } p_2(\mathcal{J}_2)\leq 0,
        \end{cases}
    \end{align}
    where
   \begin{align*}
   	p_1(\mathcal{J}_2) &= \begin{cases}
   		-p+\frac{d(d+1)}{2}+1+\frac{1-\alpha}{d-1}&\text{ if } p_2(\mathcal{J}_2)> 0,\\
   		\frac{\alpha-2}{2}p + \frac{3d}{2} + \frac{1}{(d-1)}(-\frac{\alpha p}{2} +d)  &\text{ if } p_2(\mathcal{J}_2) \leq 0.
   	\end{cases}
   \end{align*}
    If $p_2(\mathcal{J}_2) \leq 0$, then $p_1(\mathcal{J}_2) <0$ by the previous calculation when $p_2(\mathcal{J}_1) >0$. If $p_2(\mathcal{J}_2)> 0$, then
    \[
        p_1(\mathcal{J}_2) < (1-\alpha)\left(  \frac{1}{d-1} -\frac{d(d+1)}{2\alpha} \right) <0.
    \]
    Therefore, we get
    \begin{align}\label{eq_sum_J2}
         \mathcal{J}_2 \lesssim_\e M_{j}^{\e+ \frac{d}{2(d-1)}(-\alpha p +3(d-1) +2\alpha)}+M_{j}^{\e-\alpha p+\frac{d(d+1)}{2}+\alpha}.
    \end{align}
    Plugging \eqref{eq_sum_J1} and \eqref{eq_sum_J2} into \eqref{eq_v_j}, we obtain that
    \[
        \norm{\widehat{\nu}_{j+1} - \widehat{\nu}_j }_p^p\lesssim_\e M_j^{\e -\frac{\alpha p}{2}+d } + M_j^{\e -\alpha p + \frac{d(d+1)}{2}+\alpha} +M_j^{\e+ \frac{d}{2(d-1)}(-\alpha p +3(d-1) +2\alpha)}   .
    \]
    Since $p >p_\alpha$, all the exponents of $M_j$ are negative, provided that $\e>0$ is small enough. Especially, the exponent of $M_j$ in the last term is negative since \eqref{eq_rom2} holds when $ d\geq 4$. By summing over $j$, we have
    \[
        \norm{\widehat{\nu}_n }_p \lesssim \sum_{ 1 \leq j \leq n-1 } \norm{\widehat{\nu}_{j+1} - \widehat{\nu}_j }_p  +1 \lesssim_p 1.
    \]
\end{proof}
When $d=3$, \eqref{eq_rom2} does not hold. Thus, we will consider a slightly different decomposition of $\Omega_j(s_1)$. For each $1\leq 2^{s_2}\leq M_j 2^{s_1}$, we define
\begin{align*}
	\Omega_{3,j} (s_1,s_2):=\{\xi \in \Omega_j(s_1) : &\exists t \in E_{j} \ \text{such \ that}\\
	&|\langle \gamma^{(k)} (t), \xi \rangle | \leq (M_j2^{s_1}) 2^{s_2(k-1)}, \ \forall 1 \leq k \leq 3 \}.
\end{align*}
The definition of $\Omega_{3,j}(s_1,s_2) $ requires the existence of $t$ such that 
\begin{equation}\label{Mj2s1s2}
 |\langle \gamma^{(k)} (t), \xi \rangle | \leq (M_j 2^{s_1}) 2^{s_2(k-1)}   ,
\end{equation}
for all $1 \leq k \leq 3$ However,  in the definition of $\Omega_j(s_1,s_2)$, the numbers $t$ such that \eqref{Mj2s1s2} holds with $k=1$ and $t'$ such that \eqref{Mj2s1s2} holds with $k=2,3$ do not need to be the same. Hence,
$\Omega_{3,j}(s_1,s_2)$ could be strictly contained in $\Omega_j(s_1,s_2)$ for some $s_2$.

Now let us introduce the estimates for the case $d=3$.
\begin{lemma}\label{lem_mup_meas_d=3}
    For any $\e >0$, a.s., we have
    \begin{equation}\label{Omega_meas_3}
            |\Omega_{3,j}(s_1,s_2)| \lesssim_\e M_j^\e (M_j2^{s_1})^3 2^{s_1(1-\alpha)+s_2( 3 +\alpha)}  .
    \end{equation}  
\end{lemma}
Next, we define the set $\widetilde{\Omega}_{3,j+1}(s_1,s_2)$ as follows:
\begin{align*}
	\widetilde{\Omega}_{3,j+1}(s_1,s_2) &= \Omega_{3,j+1}(s_1,s_2) \backslash (\Omega_{3,j+1}(s_1,s_2-1) \cup \Omega_{j+1}(s_1-1)) \quad \text{if } \ s_1, s_2 \neq 0,\\
	\widetilde{\Omega}_{3,j+1}(s_1,0) &=\Omega_{3,j+1}(s_1,0) \backslash \Omega_{j+1}(s_1-1) \quad \text{if } \ s_1 \neq 0, \\
	\widetilde{\Omega}_{3,j+1}(0,s_2) &= \Omega_{3,j+1}(0,s_2) \backslash \Omega_{3, j+1}(0,s_2-1) \quad \text{if } \ s_2 \neq 0,\\
	\widetilde{\Omega}_{3,j+1}(0,0) &= \Omega_{3,j+1}(0,0).
\end{align*}
\begin{lemma}\label{lem_mup_osci_d=3}
    For any $\e >0$, if $\xi \in \widetilde{\Omega}_{3,j+1}(s_1,s_2)$, a.s.  we have 
    \begin{equation}\label{Omega_osci2}
        |\widehat{\nu}_{j+1}(\xi) - \widehat{\nu}_j(\xi) | \lesssim_\e (M_j 2^{s_1})^\e M_j^{-\frac{\alpha}{2}} 2^{-s_1(1-\frac{\alpha}{2})} 2^{-s_2\frac{\alpha}{2}}.
    \end{equation}
\end{lemma}
Equipped with these two lemmas, we can prove Proposition \ref{Prop_optimal} when $d=3$.

\begin{proof}[Proof of Proposition \ref{Prop_optimal} when $d = 3$] 
    Similarly, it suffices to show that for any $n\geq 1$, 
    \[
        \norm{\widehat{\nu}_n }_p \lesssim \sum_{ 1 \leq j \leq n-1 } \norm{\widehat{\nu}_{j+1} - \widehat{\nu}_j }_p  +1\lesssim 1.
    \]
    For fixed $j\geq 1$, we have
    \begin{align}\label{eq_v_j_p}
        \norm{\widehat{\nu}_{j+1} - \widehat{\nu}_j }_p^p \leq 
        \sum_{\substack{1 \leq s_1 \\ 1 \leq 2^{s_2} \leq 2^{s_1} M_j }}\int_{\widetilde{\Omega}_{3,j+1}(s_1,s_2)} |\widehat{\nu}_{j+1} - \widehat{\nu}_j|^p d\xi .
    \end{align}
    Applying Lemmas \ref{lem_mup_meas_d=3} and \ref{lem_mup_osci_d=3}, we get
    \begin{equation*}
    \begin{split}
       \int_{\widetilde{\Omega}_{3,j+1}(s_1,s_2)} |\widehat{\nu}_{j+1} - \widehat{\nu}_j|^p d\xi \lesssim_\e 
        (M_j 2^{s_1})^\e M_j^{3-\frac{\alpha p}{2}}2^{s_1(\frac{\alpha p}{2}-p+4-\alpha)}2^{s_2(3+\alpha-\frac{\alpha p}{2})}.
    \end{split}
    \end{equation*}
For fixed $s_1$, denote the exponent of $2^{s_2}$ by $p_2$, i.e., $p_2 = 3+\alpha - \frac{\alpha p}{2}$. Then, we can check that the sum over $1 \leq 2^{s_2} \leq 2^{s_1} M_j$ is bounded by
    \begin{align}\label{sum_s2_3}
            &\lesssim_\e \begin{cases}
             (M_j 2^{s_1})^\e M_{j}^{-\alpha p +6+\alpha}2^{s_1(7-p)}  &\text{ if } p_2>0,\\
            (M_j 2^{s_1})^\e M_{j}^{-\frac{\alpha p }{2} +3  }2^{s_1(\frac{\alpha p}{2}-p+4-\alpha)}&\text{ if } p_2\leq 0.
        \end{cases}
    \end{align}
    Since $p > p_\alpha$, the exponent of $2^{s_1}$ and $M_j$ in \eqref{sum_s2_3} are all negative for sufficiently small $\e >0$. Thus, we get
    \[
        \norm{\widehat{\nu}_{j+1} - \widehat{\nu}_j }_p^p \lesssim_\e M_j^\e ( M_j^{-\alpha p +6+ \alpha} + M_j^{- \frac{\alpha p}{2} +3}).
    \]
    Summing over $j$, we have the desired estimate, which is
     \[
     \begin{split}
         \norm{\widehat{\nu}_n }_p &\lesssim_\e \sum_{ 1 \leq j \leq n-1 } \norm{\widehat{\nu}_{j+1} - \widehat{\nu}_j }_p  +1\lesssim_\e 1.
     \end{split}
    \]
\end{proof}

\subsection{Estimates on the measures of sets}
Before we prove Lemmas \ref{lem_mup_meas} and \ref{lem_mup_meas_d=3}, we will introduce some necessary results.
    \begin{lemma}\label{lem_est_mu_j}
        A.s. we have $\norm{\mu_j} \lesssim 1$ uniformly in $j$.
    \end{lemma}
    Lemma \ref{lem_est_mu_j} easily follows from the Martingale convergence theorem. For example, see \cite[p. 28]{Ryou23}. 
    \begin{lemma}[Lemma 5.10 in \cite{Ryou23}]\label{lem_est_mu_ball_r}
         For any $\e>0$, a.s. we have 
        \[
        \mu_j(B(t,r)) \lesssim_\e M_j^\e r^\alpha, \qquad \forall t\in [0,1], \forall r\geq M_j^{-1}.
        \]
    \end{lemma}
    We also need an estimate of the volume of a parallelotope.
    \begin{lemma}\label{lem_vol_P}
    Assume that a set vectors $\{b_k \in \R^d : 1 \leq k \leq d\}$ is linearly independent. Let $P$ be a  parallelotope defined by
    \[
    \{ \xi \in \R^d : |\langle b_k, \xi \rangle | \leq 1/2 , 1 \leq k \leq d \}. 
    \]
    Then, $|P|= |\det(b_1,\cdots, b_d)|^{-1}$ where $\det(b_1, \cdots, b_d)$ is the determinant of a matrix whose columns are $b_1, \cdots, b_d$.
    \end{lemma}
    \begin{proof}
        Let $A$ be a matrix whose columns are $b_1, \cdots, b_d $. By the change of variable $\eta = A\xi$, we obtain that
        $$
        |P| = \int_{|\eta_i | \leq 1/2, 1 \leq i \leq d} 1 d\xi = \int_{|\eta_i|\leq 1/2} |\det(A)|^{-1}d\eta = |\det(A)|^{-1}.
        $$
    \end{proof}
    Lastly, we will use Taylor expansion of $\gamma(t)$ often. Thus, we write here. For fixed $\xi$,
    \begin{equation}\label{eq_taylor}
        \langle \gamma^{(k)} (t),\xi \rangle = 
        \sum_{n=0}^{d-k}\frac{1}{n!}\langle\gamma^{(k+n)}(s),\xi\rangle (t-s)^n, \qquad \forall t,s\in [0,1]\,, 1\leq k\leq d.
    \end{equation}

Let us explain the idea of the proof of Lemma \ref{lem_mup_meas} first. We define the parallelotope  $P(t_1,t_2)$ by the set of $\xi \in \R^d$ such that $|\langle \gamma^{(k)} (t_1), \xi \rangle | \leq 10(M_j2^{s_1})$ and $|\langle \gamma^{(k)} (t_2), \xi \rangle | \leq 10(M_j2^{s_1})$ for all $2 \leq k \leq d$.

The set $\Omega_j(s_1,s_2)$ is contained in the union of $P(t_1,t_2)$ where $t_1, t_2 \in E_j$. When $s_2$ is small, i.e. $2^{s_2} \leq (2^{s_1} M_j)^{\frac{1}{d-1}}$, estimating the measure of the union of $P(t_1,t_2)$ gives a good upper bound on $\Omega_j(s_1,s_2)$, which is \eqref{eq_Omega_meas1}. However, if $s_2$ is large, i.e. $2^{s_2} \geq (2^{s_1} M_j)^\frac{1}{d-1}$, and $|t_1-t_2|$ is large, then a small portion of $P(t_1,t_2)$ intersects $\Omega_j(s_1)$. Thus, we need an additional argument for this case.

Therefore, we will prove \eqref{eq_Omega_meas1} and \eqref{eq_Omega_meas2} in Lemma \ref{lem_mup_meas} separately. 

\begin{proof}[Proof of \eqref{eq_Omega_meas1} in Lemma \ref{lem_mup_meas}]
    Recall that our goal is to prove that a.s. we have
    \begin{equation*}
        |\Omega_j(s_1,s_2) | \lesssim_\e M_j^\e (M_j 2^{s_1})^d 2^{s_2(\frac{d(d-1)}{2}+\alpha d)} 2^{s_1\frac{d(1-\alpha)}{d-1}}.
        \end{equation*}
    We will consider three cases:
    \begin{itemize}
        \item[(a)] $2^{s_2} \leq M_j^\frac{1}{d-1}$,
        \item[(b)] $M_j^{\frac{1}{d-1}} \leq 2^{s_2} \leq \min \{M_j, (2^{s_1 }M_j)^{\frac{1}{d-1}}\}$,
        \item[(c)] $M_j \leq 2^{s_2} \leq (2^{s_1} M_j)^{\frac{1}{d-1}}$.
    \end{itemize}

    Case (a): Let $i_1, i_2$ be two integers such that
    \begin{align}\label{eq_i1_i2}
        M_{i_1-1} \leq 2^{s_2(d-1)} \leq M_{i_1}\,,\qquad M_{i_2-1} \leq 2^{s_2} \leq M_{i_2}.
    \end{align}
    For each $\xi\in \Omega_j(s_1,s_2)$, recall that there exists $t_1,t_2\in E_j$ satisfying
    \begin{align}\label{eq_property_xi_Omega}
        |\langle \gamma^{(1)}(t_1), \xi \rangle| \leq M_j2^{s_1},
    \end{align}
    and
    \begin{align}\label{eq_property_xi_Omega_2}
         \ |\langle\gamma^{(k)}(t_2), \xi \rangle| \leq M_j 2^{s_1} 2^{s_2(k-1)} \,, \qquad \forall 2 \leq k \leq d.
    \end{align}
    Let $T_1$ and $T_2$ be the sets of such points $t_1$ and $t_2$ respectively, i.e., $T_1$ is the set of $t_1 \in E_j$ such that there is $\xi \in \Omega_j(s_1,s_2)$ which satisfies \eqref{eq_property_xi_Omega} and $T_2$ is the set of $t_2 \in E_j$ such that there is $\xi \in \Omega_j(s_1,s_2)$ which satisfies \eqref{eq_property_xi_Omega_2}. 
    
    Observe that $i_1,i_2\leq j$. We will consider levels $i_1$ and $i_2$ of our set, namely $E_{i_1}$ and $E_{i_2}$, respectively. Recall that $E_{i_1}$ is the union of intervals $I_1$ of length $ M_{i_1}^{-1}$, and $E_j\subseteq E_{i_1}$. Denote $\widetilde{T}_1$ by the set of the left endpoints of intervals $I_{1}$ which intersect $T_1$. Similarly, $E_{i_2}$ is the union of intervals $I_{2}$ of length $ M_{i_2}^{-1}$ with $E_j\subseteq E_{i_2}$. Denote $\widetilde{T}_2$ by the set of the left endpoints of intervals $I_2 $ which intersect $T_2$.

    For each $\tilde{t}_1\in \widetilde{T}_1$ and $ \tilde{t}_2 \in \widetilde{T}_2$, we define
    \[
    \begin{split}
    \Omega_j(s_1,s_2,\tilde{t}_1 , \tilde{t}_2)=\{ \xi \in \R^d : &|\langle \gamma^{(1)}(\tilde{t}_1), \xi \rangle| \leq 100dM_j2^{s_1},\\
    &|\langle \gamma^{(k)}(\tilde{t}_2), \xi \rangle| \leq 100dM_j 2^{s_1} 2^{s_2(k-1)} \,, \forall 2 \leq k \leq d \}.
    \end{split}
    \]
    We claim that
    \begin{align}\label{eq_union_Omegajs1s2}
        \Omega_j(s_1,s_2) \subseteq \bigcup\limits_{\tilde{t}_1 \in \widetilde{T}_1, \tilde{t}_2 \in \widetilde{T}_2} \Omega_j (s_1,s_2, \tilde{t}_1, \tilde{t}_2)\,.
    \end{align}
    
    Indeed, in each interval $I_1$, we have that $|\tilde{t}_1 -t_1| \leq M_{i_1}^{-1} \leq 2^{-s_2(d-1)} $ for some $\xi \in \Omega_j(s_1,s_2)$. By Taylor expansion \eqref{eq_taylor} with $t=\tilde{t}_1$ and $s=t_1(\xi)$, we get 
    \[
    |\langle \gamma^{(1)}(\tilde{t}_1), \xi \rangle | \leq dM_j 2^{s_1}.
    \]
    Similarly, using the fact that $|\tilde{t}_2 -t_2| \leq M_{i_2}^{-1} \leq 2^{-s_2}$ and \eqref{eq_taylor}, we have 
    \[
    |\langle\gamma^{(k)}(\tilde{t}_2), \xi \rangle | \leq dM_j2^{s_1+s_2(k-1)},\qquad \forall2 \leq k \leq d.
    \]
    Thus, we have \eqref{eq_union_Omegajs1s2}. Now, note that
    \[
    \det(\gamma^{(1)}(\tilde{t}_1), \gamma^{(2)}(\tilde{t}_2), \cdots ,\gamma^{(d)}(\tilde{t}_2)) \gtrsim 1
    \]
    uniformly in $\tilde{t}_1\in \widetilde{T}_1$ and $\tilde{t}_2\in \widetilde{T}_2$. Thus it follows from Lemma \ref{lem_vol_P} that
    \begin{align}\label{eq_est_Omega_meas_a}
        |\Omega_j(s_1,s_2,\tilde{t}_1, \tilde{t}_2)| \lesssim (M_j 2^{s_1})^d 2^{s_2 \frac{d(d-1)}{2}}.
    \end{align}   
    Since $M_{i_1} \lesssim_\e M_{i_1-1}^{1+\e} \lesssim 2^{s_2(d-1)(1+\e)}$ and, $M_{i_2} \lesssim_\e M_{i_2-1}^{1+\e} \lesssim 2^{s_2(1+\e)} $, using \eqref{cond_b_jM_j} and Lemma \ref{lem_est_mu_j}, a.s. we obtain that
    \begin{align}\label{eq_est_T1_a}
        \# \widetilde{T}_1 \lesssim \norm{\mu_{i_1}} \beta_{i_1}^{-1}M_{i_1} \lesssim_\e  2^{s_2(d-1)\alpha +\e},
    \end{align}
    and
    \begin{align}\label{eq_est_T2_a}
        \# \widetilde{T}_2 \lesssim \norm{\mu_{i_2}} \beta_{i_2}^{-1}M_{i_2} \lesssim_\e  2^{s_2\alpha +\e}.
    \end{align}
    Combining the above estimates, we obtain that a.s. \eqref{eq_Omega_meas1} holds.

    Case (b): We can prove \eqref{eq_Omega_meas1} in a similar way with some modifications. Since $2^{s_2} \leq M_j \leq 2^{s_2(d-1)}$, there exist integers $i_1, i_2$ satisfying \eqref{eq_i1_i2} with $i_2\leq j\leq i_1$.
    We consider the same sets $E_{i_2}$ and $\widetilde{T}_2$ as in case $(a)$. However, we cannot use the same $\widetilde{T}_1$ in case (a). Since $E_{i_1} \subseteq E_j$, there could exist $t_1 \in T_1$ in $ E_j$ but not in $E_{i_1}$.
    
    Therefore, we divide $E_j$ into intervals $I_{1}$ of length $M_{i_1}^{-1}$, and we use this new $I_1$ to define $\widetilde{T}_1$, i.e., $\widetilde{T}_1$ is the set of left endpoints of the intervals $I_1 \subset E_j $ of length $M_{i_1}^{-1}$  which intersect $T_1$. Then, a.s. we get
   \begin{equation}\label{eq_est_T1_b}
    \#\widetilde{T}_1 \lesssim \norm{\mu_j} \beta_j^{-1} M_{i_1} \lesssim_\e M_j^{\alpha-1+\e } 2^{s_2(d-1)}.     
    \end{equation}
    Combining \eqref{eq_est_Omega_meas_a}, \eqref{eq_est_T2_a} and \eqref{eq_est_T1_b}, a.s. we have
    $$
    |\Omega_j(s_1,s_2)| \lesssim_\e M_j^\e (M_j 2^{s_1})^d 2^{s_2 \frac{d(d-1)}{2}} M_j^{\alpha-1} 2^{s_2(d-1+\alpha)}.
    $$
    Due to the range of $s_2$, we have $M_j^{-1} 2^{s_2(d-1)} \leq 2^{s_1} \leq 2^{s_1 \frac{d}{d-1}}$ and it implies \eqref{eq_Omega_meas1}.
    
    Case (c): Since $M_j \leq 2^{s_2} $ and $M_j \leq 2^{s_2(d-1)}$, there exist $i_1,i_2\geq j$ such that \eqref{eq_i1_i2} holds. Instead of $E_{i_1}$ and $E_{i_2}$ in case (a), we divide $E_j$ into intervals $I_1$ of length $M_{i_1}^{-1}$, and into intervals $I_2$ of length $M_{i_2}^{-1}$ respectively. The sets $\widetilde{T}_1$ and $\widetilde{T}_2$ are then defined in the same way in $(a)$ by using these new $I_1$ and $I_2$. These decompositions a.s. give us \eqref{eq_est_Omega_meas_a}, \eqref{eq_est_T1_b} and 
    \[
    \#\widetilde{T}_2 \lesssim \norm{\mu_j} \beta_j^{-1} M_{i_2}^{-1} \lesssim_\e M_j^{\alpha-1+\e } 2^{s_2}.
    \]
    Thus, a.s. we have 
    $$
    |\Omega_j(s_1,s_2)| \lesssim_\e M_j^\e (M_j 2^{s_1})^d 2^{s_2 \frac{d(d-1)}{2}} M_j^{-2(1-\alpha)} 2^{s_2d}.
    $$
    Due to the range of $s_2$, we have $M_j^{-2}2^{s_2d} \leq 2^{s_1 \frac{d}{d-1}} M_j^{-\frac{d-2}{d-1}} \leq 2^{s_1 \frac{d}{d-1}}$. Thus, we conclude that \eqref{eq_Omega_meas1} holds a.s.
\end{proof}
Now we turn to the proof of \eqref{eq_Omega_meas2}.
\begin{proof}[Proof of \eqref{eq_Omega_meas2} in Lemma \ref{lem_mup_meas}]
Recall that our goal is to prove that a.s. we have
\begin{equation*}
        |\Omega_j(s_1,s_2)| \lesssim_\e M_j^\e (M_j 2^{s_1})^d 2^{s_2 \frac{d(d-1)}{2} } M_j^\alpha 2^{s_1} ( M_j 2^{s_1-s_2})^{\frac{\alpha}{d-2}}2^{s_1\frac{1-\alpha}{d-1}}.
    \end{equation*}
We will consider three cases:
\begin{itemize}
    \item[(a)] $(2^{s_1} M_j)^{\frac{1}{d-1}} \leq 2^{s_2} \leq  M_j$,
    \item[(b)] $M_j \leq 2^{s_2} \leq 2^{s_1/(d-1)} M_j$,
    \item[(c)] $ 2^{s_1/(d-1)} M_j \leq 2^{s_2} \leq 2^{s_1} M_j$.
\end{itemize}

    Case (a): For each $\xi\in \Omega_j(s_1,s_2)$, recall that there exists $t_1,t_2\in E_j$ satisfying
    $$
    |\langle \gamma^{(1)}(t_1), \xi \rangle| \leq M_j2^{s_1}\,,\qquad |\langle \gamma^{(2)}(t_1), \xi \rangle| \leq (M_j2^{s_1})^2
    $$
    and
    $$
    |\langle\gamma^{(k)}(t_2), \xi \rangle| \leq M_j 2^{s_1} 2^{s_2(k-1)} \,, \forall 2 \leq k \leq d.
    $$
    Let $T_1$ and $T_2$ be the sets of such points $t_1$ and $t_2$ respectively. By the assumption on $s_2$, there is an integer $i$ such that $M_{i-1} \leq 2^{s_2} \leq M_i$. 
    Decompose $E_{j}$ into intervals $I_{1}$ of length $ (M_j 2^{s_1})^{-1}$, and define $\widetilde{T}_1$ as in the previous argument, i.e., $\widetilde{T}_1$ is the set of left endpoints of intervals $I_1$ which intersect $T_1$.
    Next, we consider level $i$ of our set, namely $E_i$, which is the union of intervals $I_2$ of length $M_i^{-1}$, and we define the corresponding set $\widetilde{T}_2$ in the same way.
     
    For each $\tilde{t}_1\in \widetilde{T}_1$ and $ \tilde{t}_2 \in \widetilde{T}_2$, we define
    \begin{align*}
        \Omega_j(s_1,s_2,\tilde{t}_1 , \tilde{t}_2)=\{x \in \R^d : &|\langle \gamma^{(1)}(\tilde{t}_1), \xi \rangle| \leq 100dM_j2^{s_1}\,, \\
        &|\langle \gamma^{(2)}(\tilde{t}_1), \xi \rangle| \leq 100d(M_j2^{s_1})^2\,,\\ 
        &|\langle\gamma^{(k)}(\tilde{t}_2), \xi \rangle| \leq 100dM_j 2^{s_1} 2^{s_2(k-1)} \,,\quad \forall 2 \leq k \leq d \}.
    \end{align*}
    Since $|I_1| = (M_j2^{s_1})^{-1}$ and $|I_2| = M_i^{-1} \leq 2^{-s_2}$, we apply Taylor expansion \eqref{eq_taylor} and we can check that
    \begin{align*}
        \Omega_j(s_1,s_2) \subseteq \bigcup\limits_{\tilde{t}_1 \in \widetilde{T}_1, \tilde{t}_2 \in \widetilde{T}_2} \Omega_j (s_1,s_2, \tilde{t}_1, \tilde{t}_2).
    \end{align*}
    To estimate the measure of the set on the right-hand side, we decompose it further. Let $A_0$ be a union of $\Omega_j(s_1,s_2,\tilde{t}_1, \tilde{t}_2)$ such that $|\tilde{t}_1 - \tilde{t}_2| \lesssim 2^{-s_2} ( M_j 2^{s_1-s_2} )^{\frac{1}{d-2}}$. For $l \geq 1$, let $A_l$ be a union of $\Omega_j(s_1,s_2,\tilde{t}_1, \tilde{t}_2)$ such that $|\tilde{t}_1 - \tilde{t}_2| \sim 2^{l-s_2} ( M_j 2^{s_1-s_2} )^{\frac{1}{d-2}}$. Therefore, we have 
    $$\Omega_j (s_1,s_2) \subseteq \bigcup_{l \geq 0 } A_l.$$
    
    Let us first estimate $|A_0|$. Since
    $$
    \det(\gamma^{(1)}(\tilde{t}_1), \gamma^{(2)}(\tilde{t}_2), \cdots ,\gamma^{(d)}(\tilde{t}_2)) \gtrsim 1
    $$
    uniformly in $\tilde{t}_1$ and $\tilde{t}_2$, Lemma \ref{lem_vol_P} implies that 
    \[
    |\Omega_j(s_1,s_2,\tilde{t}_1, \tilde{t}_2)| \lesssim (M_j 2^{s_1})^d 2^{s_2 \frac{d(d-1)}{2}},\qquad \forall \tilde{t}_1\in \widetilde{T}_1,\tilde{t}_2\in \widetilde{T}_2.    \]
    Here, we did not use the inequality $|\langle \gamma^{(2)}(\tilde{t}_1), \xi \rangle| \leq 100d(M_j2^{s_1})^2$. We will use it to estimate $|A_l|$ for $l \geq 1$.
    
    For fixed $\tilde{t}_1\in \widetilde{T}_1$, since $M_i^{-1} \leq 2^{-s_2}$ and $M_j^{-1} \leq 2^{-s_2} (M_j 2^{s_1-s_2})^{\frac{1}{d-2}}$, we apply Lemma \ref{lem_est_mu_ball_r}. Then, a.s. we have
    \[
    \begin{split}
    \#\{\tilde{t}_2 \in \widetilde{T}_2: |\tilde{t}_1 - \tilde{t}_2| \lesssim 2^{-s_2} (M_j 2^{s_1-s_2})^\frac{1}{d-2}\} &\lesssim \beta_i^{-1} M_i \mu_i \big(B(\tilde{t}_1, 2^{-s_2} ( M_j 2^{s_1-s_2})^\frac{1}{d-2})\big)\\
    &\lesssim_\e M_i^{\alpha+\e } \big(2^{-s_2} (M_j 2^{s_1-s_2})^\frac{1}{d-2}\big))^\alpha\\
    &\lesssim_\e M_i^{\e} (M_j 2^{s_1-s_2})^\frac{\alpha}{d-2}.
    \end{split}
    \]
    Additionally, a.s. we have
    \[
    \#\widetilde{T}_1 \lesssim \norm{\mu_{j}} \beta_{j}^{-1}M_{j} 2^{s_1} \lesssim_\e  M_j^{\alpha +\e }2^{s_1}.
    \]
    Therefore, we obtain that
    \begin{equation}\label{Omega_meas_A1}
        |A_1|\lesssim_\e M_j^\e (M_j 2^{s_1})^d 2^{s_2 \frac{d(d-1)}{2}} M_j^\alpha 2^{s_1} (M_j 2^{s_1-s_2})^\frac{\alpha}{d-2}.
    \end{equation}
    Next, we estimate $|A_l|$. Observe that $\Omega_j(s_1,s_2,\tilde{t}_1, \tilde{t}_2)$ is contained in the  parallelotope
    \[
    \begin{split}
    \{ \xi \in \R^d : &|\langle \gamma^{(1)}(\tilde{t}_1), \xi \rangle| \leq 100dM_j2^{s_1}, |\langle \gamma^{(2)}(\tilde{t}_1)\,, \xi \rangle| \leq 100d(M_j2^{s_1})^2,\\
    &\hspace{2cm}|\langle\gamma^{(k)}(\tilde{t}_2), \xi \rangle| \leq 100dM_j 2^{s_1} 2^{s_2(k-1)} \,, \forall 2 \leq k \leq d-1 \}.
    \end{split}
    \]
    By Lemma \ref{lem_vol_P}, we can bound $|\Omega_j(s_1,s_2, \tilde{t}_1 , \tilde{t}_2)| $ by a constant multiple of 
    \begin{equation}  \label{eq_Omega_meas_A2_est}
    \begin{split}
    (M_j 2^{s_1})^{d+1} 2^{s_2 \frac{(d-2)(d-1)}{2}}|\det (\gamma^{(1)} (\tilde{t}_1),\gamma^{(2)} (\tilde{t}_1), \gamma^{(2)} (\tilde{t}_2), \cdots, \gamma^{(d-1)} (\widetilde{t_2}))|^{-1}.
    \end{split}
    \end{equation}
    By the Taylor expansion \eqref{eq_taylor}, we have
    \[
    \gamma^{(2)} (\tilde{t}_1) =  \gamma^{(2)} (\tilde{t}_2 )  + \gamma^{(3)} (\tilde{t}_2 ) (\tilde{t}_1 - \tilde{t}_2 ) + \cdots +\frac{1}{(d-2)!} \gamma^{(d)} (\tilde{t}_2 )  (\tilde{t}_1 - \tilde{t}_2 )^{d-2}.
    \]
    Thus, the determinant in \eqref{eq_Omega_meas_A2_est} is bounded by a constant multiple of
    \[
    \begin{split}
        |\det (\gamma^{(1)} (\tilde{t}_1), \gamma^{(2)} (\tilde{t}_2), \cdots, \gamma^{(d-1)} (\tilde{t}_2), \gamma^{(d)} (\tilde{t}_2))| \cdot |\tilde{t}_1 -\tilde{t}_2 |^{d-2}\sim |\tilde{t}_1 -\tilde{t}_2|^{d-2}.
    \end{split}
    \]
    Plug this into \eqref{eq_Omega_meas_A2_est}, we get
    \[
    \begin{split}
    |\Omega_j(s_1,s_2, \tilde{t}_1 , \tilde{t}_2)| \lesssim &(M_j 2^{s_1})^{d} 2^{s_2 \frac{d(d-1)}{2}} \cdot 2^{-l(d-2) }  .
    \end{split}
    \]
    By the same calculation as in the case of $A_0$, a.s. we have that $\#\widetilde{T}_1 \lesssim_\e  M_j^{\alpha +\e }2^{s_1}$, and for any $\tilde{t}_1\in \widetilde{T}_1$,
    \begin{align}\label{eq_card_T2_a}
        \#\{\tilde{t}_2 \in \widetilde{T}_2: |\tilde{t}_1 - \tilde{t}_2 | \lesssim 2^{l-s_2} (M_j 2^{s_1-s_2})^\frac{1}{d-2}\} \lesssim_\e M_i^{\e} 2^{l\alpha}(M_j 2^{s_1-s_2})^\frac{\alpha}{d-2}.
    \end{align}
    Combining the estimates above, a.s. we have
    \begin{equation}\label{Omega_meas_A2}
    \begin{split}
        |A_{l}| \lesssim_\e 2^{-l(d-2-\alpha)} M_j^\e (M_j 2^{s_1})^d 2^{s_2 \frac{d(d-1)}{2}} M_j^\alpha 2^{s_1} (M_j 2^{s_1-s_2})^\frac{\alpha}{d-2},\qquad \forall l\geq 1.
    \end{split}
    \end{equation}
    Since $d-2-\alpha>0$, using \eqref{Omega_meas_A1} and \eqref{Omega_meas_A2}, we obtain that
    \begin{align}
        \left|\bigcup_{l \geq 0 } A_l\right| \leq \sum_{l\geq 0}|A_{l}|\lesssim_\e M_j^\e (M_j 2^{s_1})^d 2^{s_2 \frac{d(d-1)}{2}} M_j^\alpha 2^{s_1} (M_j 2^{s_1-s_2})^\frac{\alpha}{d-2}.
    \end{align}
    Therefore, we have established \eqref{eq_Omega_meas2}.
    
    Case $(b)$: By the range of $s_2$, we have $2^{-s_2} \leq M_j^{-1} \leq 2^{-s_2} ( M_j 2^{s_1-s_2})^\frac{1}{d-2}$. Instead of $E_i$ in case (a), we decompose $E_j$ into intervals $I_2$ of length $ \sim 2^{-s_2}$ and define $\widetilde{T}_2$ in the same way. Then, we can follow the previous argument in case (a) line by line except that the estimate \eqref{eq_card_T2_a} is replaced by 
    \[
    \#\{\tilde{t}_2 \in \widetilde{T}_2: |\tilde{t}_1 - \tilde{t}_2 | \lesssim 2^{l-s_2} (M_j 2^{s_1-s_2})^\frac{1}{d-2}\}
    \lesssim_\e M_i^{\e } 2^{l\alpha}  (M_j 2^{s_1-s_2})^\frac{\alpha}{d-2}(2^{s_2 } M_j^{-1})^{1-\alpha} ,
    \]
    for any $l \geq 0$. Since $2^{s_2} M_j^{-1} \leq 2^{s_1/(d-1)}$, a.s. we have \eqref{eq_Omega_meas2}.

    Case $(c)$: Since $2^{-s_2} \leq M_j^{-1} $, as in case (b), we divide $E_j$ into intervals $I_2$ of length $\sim 2^{-s_2}$ and define $\widetilde{T}_2$ in the same way. Again, we can follow the same argument for case (a) except that \eqref{eq_card_T2_a} is replaced by
    \[
    \#\{\tilde{t}_2 \in \widetilde{T}_2: |\tilde{t}_1 - \tilde{t}_2 | \lesssim 2^{l-s_2} (M_j 2^{s_1-s_2})^\frac{1}{d-2}\}\lesssim 2^l(M_j 2^{s_1-s_2})^\frac{1}{d-2},\qquad \forall l\geq 0.
    \]
    Here, we used that $2^{-s_2} ( M_j 2^{s_1-s_2})^\frac{1}{d-2 } \leq M_j^{-1}$ since $ 2^{s_1/(d-1)} M_j \leq 2^{s_2} $. 
    We use the fact that  $(M_j 2^{s_1-s_2})^{\frac{1}{d-2}} \leq 2^{s_1 \frac{1}{d-2}}$. Then, a.s. \eqref{eq_Omega_meas2} holds.
    \end{proof}

Let us turn to the proof of Lemma \ref{lem_mup_meas_d=3}. 
\begin{proof}[Proof of Lemma \ref{lem_mup_meas_d=3}]
Recall that we want to show that a.s. we have
\begin{equation*}
            |\Omega_{3,j}(s_1,s_2)| \lesssim_\e M_j^\e (M_j2^{s_1})^3 2^{s_1(1-\alpha)+s_2( 3 +\alpha)}.
    \end{equation*}  
We will consider two cases:
\begin{itemize}
    \item[(a)] $1\leq 2^{s_2}\leq M_j$,
    \item[(b)] $M_j\leq 2^{s_2}\leq 2^{s_1}M_j$.
\end{itemize}
Case (a):  Let $i$ be an integer such that $M_{i-1} \leq 2^{s_2} \leq M_i$.
We consider level $i$ of our set, namely $E_i$, which is the union of intervals $I_{i}$ of length $ M_{i}^{-1}$. For each $\xi\in \Omega_{3,j}(s_1,s_2)$, there exists $t \in E_j$ such that
        \[
        |\langle \gamma^{(k) }(t), \xi \rangle | \leq M_j 2^{s_1} 2^{s_2(k-1)}, \qquad \forall 1 \leq k \leq 3.
        \]
        Let $T$ be the set of such points $t$ and let $\widetilde{T}$ be the set of left endpoints of intervals $I_{i}$ which intersect $T$. For each $\tilde{t} \in \widetilde{T}$, define
        \[
        \Omega_{3,j}(s_1, s_2, \tilde{t}) = \{ \xi \in \R^d : |\langle \gamma^{(k)}(\tilde{t}), \xi \rangle | \leq 10 d(M_j 2^{s_1}) 2^{s_2(k-1)} \ \forall 1 \leq k \leq 3 \}.
        \]
        We can check that 
        \[
        \widetilde{\Omega}_{3,j} (s_1,s_2) \subseteq \bigcup\limits_{\tilde{t} \in \widetilde{T}} \Omega_{3,j}(s_1,s_2,\tilde{t}).
        \]
        By Lemma \ref{lem_vol_P}, we have
        \[
        |\Omega_{3,j} (s_1,s_2, \tilde{t} )| \lesssim (M_j 2^{(s_1+s_2)})^3.
        \]
        Also note that $\#\widetilde{T} \lesssim\norm{\mu_i} \beta_i^{-1} M_i$, since $2^{s_2} \leq M_j$.
        
        Combining the above estimates with Lemma \ref{lem_est_mu_j} and \eqref{cond_b_jM_j}, a.s. we get
        \begin{equation*}
        \begin{split}
            |\Omega_{j}(s_1,s_2)| &\lesssim (M_j 2^{(s_1+s_2)})^3\norm{\mu_i} \beta_{i}^{-1} M_i \\
            &\lesssim_\e M_j^\e (M_j2^{s_1})^3 2^{s_2(3+\alpha)}. 
        \end{split}
        \end{equation*}
    Therefore, we obtain \eqref{Omega_meas_3}.

    Case (b): Since $M_j \leq 2^{s_2} $, we divide $E_j$ into intervals of length $\sim 2^{-s_2}$ and define the set $\widetilde{T}$ in the same way. Note that $\#\widetilde{T}\lesssim \norm{\mu_j} \beta_j^{-1} 2^{s_2}$. Following the same argument in case (a), a.s. we obtain \eqref{Omega_meas_3}.
    \end{proof}

\subsection{Estimates on the oscillatory integrals}
In this section, we prove Lemma \ref{lem_mup_osci} and \ref{lem_mup_osci_d=3}. To begin with, we introduce Hoeffding's inequality.
    \begin{lemma}
         [Hoeffding's inequality \cite{Ho63}]
    Let $X_1, \cdots, X_n$ be independent real random variables such that $a_i \leq X_i \leq b_i$ and $S_n := X_1 + \cdots + X_n$. For $t>0$,
    \begin{equation}\label{eq_hoeff}
        \mathbb{P} \left( \big|S_n - \mathbb{E}(S_n) \big|>t \right) \leq 2 \exp{\left(\frac{-2t^2}{\sum_{i=1}^n (b_i-a_i)^2} \right)}.
    \end{equation}
     \end{lemma}
    For an interval $I \subseteq [0,1]$, we write
        $$X_{I} = \int_{I} \beta_{j+1} \mathbf{1}_{E_{j+1}} e(-\langle \gamma(t), \xi \rangle ) dt.$$
        Recall that $\cI_j$ is the collection of invervals of length $M_j^{-1}$ which cover $[0,1]$. We can check that
        \[
       \widehat{\nu}_{j+1}(\xi)= \sum_{I_j \in \mathcal{I}_j} X_{I_j} \qquad \text{and}\qquad  \widehat{\nu}_j (\xi)= \mathbb{E}\left(\sum_{I_j \in \cI_j} X_{I_j} | E_j \right),
        \]
        which implies that
        \begin{align*}
            |\widehat{\nu}_{j+1}(\xi) - \widehat{\nu}_j(\xi) |=\left|\sum_{I_j \in \cI_j} X_{I_j} - \mathbb{E}\left(\sum_{I_j \in \cI_j} X_{I_j} |E_j \right) \right|.
        \end{align*}
        Thus, we need an upper bound of $\sum_{I_j \in \mathcal{I}_j} |X_{I_j}|^2$ to use Hoeffding's inequality.
        \begin{lemma} \label{claim_1osc} 
            If $d\geq 4$, for each $\xi \in \widetilde{\Omega}_{j+1}(s_1,s_2)$ and $E_j$, there exists a quantity $C(E_j,s_1,s_2)$ such that 
            \begin{equation}\label{sum_X^2}
            \sum_{I_j\in\cI_j}|X_{I_j}|^2 \leq C(E_j,s_1,s_2)    
            \end{equation}
            conditioned on $E_j$ and a.s. we have
    \begin{align}\label{eq_sum_Xj}
             C(E_j,s_1,s_2) \lesssim_\e M_j^{-\alpha+\e} 2^{-s_1(2-\alpha)} 2^{-s_2\alpha}
        \end{align}
        uniformly in $\xi$ and $E_j$.
        \end{lemma}

The proof of Lemma \ref{lem_mup_osci_d=3} is the same as the proof of Lemma \ref{lem_mup_osci} except that we use the following lemma instead of Lemma \ref{claim_1osc}.

\begin{lemma} \label{claim_2osc} 
            If $d= 3$, for each $\xi \in \widetilde{\Omega}_{3,j+1}(s_1,s_2)$, there exists a quantity $C(E_j,s_1,s_2)$ such that 
            $$
            \sum_{I_j \in \mathcal{I}_j} |X_{I_j}|^2 \leq C(E_j,s_1,s_2)
            $$
            conditioned on $E_j$ and a.s. we have
    \begin{align}\label{eq_sum_Xj_d=3}
            C(E_j,s_1,s_2)  \lesssim_\e M_j^{-\alpha+\e} 2^{-s_1(2-\alpha)} 2^{-s_2\alpha}
        \end{align}
        uniformly in $\xi$ and $E_j$.
        \end{lemma}

\begin{proof}[Proof of Lemmas \ref{lem_mup_osci} and \ref{lem_mup_osci_d=3} using Lemmas \ref{claim_1osc} and \ref{claim_2osc}]
First, we prove Lemma \ref{lem_mup_osci}. We claim that it is enough to show that \eqref{Omega_osci}  holds for all $\xi \in \widetilde{\Omega}_{j+1}^\#(s_1,s_2)$, where
\begin{align*}
    \widetilde{\Omega}_{j+1}^\#(s_1,s_2)= (C_d2^{s_1}M_{j+1}^2)^{-1} \mathbb{Z}^d \cap \widetilde{\Omega}_{j+1}(s_1,s_2),
\end{align*}
for sufficiently large constant $C_d>0$. Indeed, let us assume that \eqref{Omega_osci} holds for all $\xi \in \widetilde{\Omega}^\#_{j+1}(s_1,s_2)$, namely,
    \begin{align*}
        |\widehat{\nu}_{j+1}(\xi) - \widehat{\nu}_j(\xi) |\lesssim_\e M_j^{-\frac{\alpha}{2}+\e}2^{-s_1(1-\frac{\alpha}{2})-s_2\frac{\alpha}{2}}, \qquad \forall \xi \in \widetilde{\Omega}^\#_{j+1}(s_1,s_2).
    \end{align*}
    For $h\in \R^d$, observe that
    \[
    \begin{split}
    |\widehat{\nu}_{j+1} (\xi +h) -\widehat{\nu}_{j+1}(\xi)| &=  \int \beta_{j+1} \mathbf{1}_{E_{j+1}} ( e(-\langle \gamma (t), \xi+h \rangle )-e(-\langle \gamma (t), \xi \rangle )) dt\\
    &\lesssim \beta_{j+1} |h| \lesssim M_{j+1} |h|.
    \end{split}
    \]
    Similarly, we have $|\widehat{\nu}_{j} (\xi +h) -\widehat{\nu}_{j}(\xi)|\lesssim M_{j}|h|$.
    
    If $|h| \lesssim 2^{-s_1}M_{j+1}^{-2}$, then we have
      \begin{align*}
          |\widehat{\nu}_{j+1}(\xi+h) - \widehat{\nu}_{j} (\xi +h)| 
        &\lesssim |\widehat{\nu}_{j+1}(\xi) - \widehat{\nu}_j(\xi) | + 2^{-s_1}M_j^{-1}\\
        &\lesssim_\e M_j^{-\frac{\alpha}{2}+\e}2^{-s_1(1-\frac{\alpha}{2})-s_2\frac{\alpha}{2}}+2^{-s_1}M_j^{-1}\\
        &\lesssim_\e M_j^{-\frac{\alpha}{2}+\e}2^{-s_1(1-\frac{\alpha}{2})-s_2\frac{\alpha}{2}}.
      \end{align*}
        This means that that \eqref{Omega_osci} holds for all $\xi \in \widetilde{\Omega}_{j+1}(s_1,s_2)$.
        
    Now assume that $\xi \in \widetilde{\Omega}_{j+1}^\#(s_1,s_2)$. The inequality \eqref{sum_X^2} in Lemma \ref{claim_1osc} guarantees the existence of $a_{I_j} $ and $b_{I_j}$ such that $a_{I_j} \leq X_{I_j}\leq b_{I_j}$ and $\sum_{I_j \in \mathcal{I}_j} (b_{I_j}-a_{I_j})^2 \leq 4C(E_j,s_1,s_2)$. Thus, we can use $4C(E_j,s_1,s_2)$ instead of $\sum_{I_j \in \mathcal{I}_j} (b_{I_j}-a_{I_j})^2$ in Hoeffding's inequality.
    
    Using Hoeffding's inequality with \eqref{sum_X^2} in Lemma \ref{claim_1osc}, we obtain that
    \begin{equation}\label{eq_prob_1}
        \mathbb{P} \left( |\widehat{\nu}_{j+1}(\xi) - \widehat{\nu}_j(\xi) |>C(E_j,s_1,s_2)^{1/2} (2^{s_1}M_j)^\e \big|E_j \right)
        \lesssim \exp(-(2^{s_1}M_j^{\e})^2/2).
    \end{equation}
    Note that $\widetilde{\Omega}_{j+1}(s_1,s_2) \subseteq  \Omega_{j+1}(s_1)$ and $ \Omega_{j+1}(s_1)$ is contained in a ball of radius $\lesssim(M_{j+1} 2^{s_1})^{d}$ and the elements in $\widetilde{\Omega}_{j+1}^\#(s_1,s_2)$ are $\sim 2^{-s_1} M_{j+1}^{-2}$ separated. Thus,
    \begin{align*}
        \# \widetilde{\Omega}_{j+1}^\#(s_1,s_2) \lesssim M_{j+1}^{d(d+2)} 2^{s_1d(d+1)}.
    \end{align*}
     Since \eqref{eq_prob_1} holds uniformly in $E_j$ and $\xi \in \widetilde{\Omega}_{j+1}^\#(s_1,s_2)$, we have
     \[
        \begin{split}
            \mathbb{P} \bigg( |\widehat{\nu}_{j+1}(\xi) - \widehat{\nu}_j(\xi) |&>C(E_j,s_1,s_2)^{1/2}(2^{s_1}M_j)^\e, \text{ for some } \xi \in\widetilde{\Omega}_{j+1}^\#(s_1,s_2) \bigg)\\
            &\hspace{3cm}\lesssim M_{j+1}^{d(d+2)}2^{s_1d(d+1)}\exp(-2 (2^{s_1}M_j)^{2\e}).
        \end{split}
        \]
        Since $2^{s_2} \leq 2^{s_1} M_{j+1}$, we sum both sides over $j, s_1$, and $s_2$ and we can check that the sum of the right-hand side is bounded. By Borel-Cantelli Lemma, a.s. we have
    \begin{equation}\label{eq_prob_2}
    |\widehat{\nu}_{j+1}(\xi) - \widehat{\nu}_j(\xi) | \leq C(E_j,s_1,s_2)^{1/2} (2^{s_1}M_j)^{\e}, \qquad \forall \xi\in \widetilde{\Omega}_{j+1}(s_1,s_2).    
    \end{equation}
    Combining with \eqref{eq_sum_Xj} in Lemma \ref{claim_1osc}, we have established \eqref{Omega_osci}.

    Next, let us turn to the proof of Lemma \ref{lem_mup_osci_d=3}. The proof is same with the previous argument except that $\widetilde{\Omega}_{j+1}\#(s_1,s_2)$ is replaced by 
    $$\widetilde{\Omega}_{3,j+1}^\#(s_1,s_2):=(C_d2^{s_1}M_{j+1}^2)^{-1})\cap \widetilde{\Omega}_{3,j+1}(s_1,s_2)$$
    and that we use Lemma \ref{claim_2osc} instead of Lemma \ref{claim_1osc}.
\end{proof}

To establish Lemma \ref{claim_1osc} and \ref{claim_2osc}, we use the following lemma on the oscillatory integrals.
    \begin{lemma}[Theorem 1.1 in \cite{ACK_book}]\label{lem_osci_est}
        Suppose that $d \geq 1$, and $I$ is an interval in $[0,1]$. Then,
        \[
        \left|\int_I e(-\langle \gamma(t), \xi \rangle ) dt \right| \lesssim \min \{ |I|, H_I(\xi)^{-1}\} , \qquad \forall \xi\in \R^d,
        \]
        where
        \[
        H_I(\xi) = \inf_{t \in I} \sum_{k=1}^d \bigg|\frac{1}{k!}\langle \gamma^{(k)}(t),\xi \rangle \bigg|^{1/k}.
        \]
    \end{lemma}
    The following estimate holds regardless of whether $d\geq 4$ or $d = 3$.
    \begin{lemma}\label{lem_XI_s_2=0}
        For each $\xi \in \widetilde{\Omega}_{j+1}(s_1,s_2)$, we have
        \begin{align*}
            |X_{I_j}|  \lesssim \beta_{j+1} M_j^{-1} 2^{-s_1}.
        \end{align*}
    \end{lemma}
    \begin{proof}
        If $s_1=0$, we use that the length of $I_j$ is $M_j^{-1}$. Thus, we obtain that 
        \begin{align*}\label{eq_X_estimate_1}
             |X_{I_j} | \lesssim \beta_{j+1} M_j^{-1}. 
        \end{align*}
        If $s_1 \neq 0$, $\xi \in \widetilde{\Omega}_{j+1}(s_1,0)$ implies that $\xi \notin \Omega_{j+1}(s_1-1)$. Thus, for any $t \in E_{j+1}$, there exists $1 \leq k \leq d$ such that
         \[
         |\langle \gamma^{(k)} (t), \xi \rangle | > (M_{j+1}2^{s_1})^k.
         \]
         This yields that for each interval $I_{j+1} \in \mathcal{I}_{j+1}$ contained in $ E_{j+1}$, we have
         \begin{align*}\label{eq_Hxi}
             H_{I_{j+1}}(\xi)\gtrsim M_{j+1}2^{s_1}.
         \end{align*}
         For each $I_j \in \mathcal{I}_j$, at most $m_{j+1}$ such intervals are contained in $I_j \cap E_{j+1}$. Thus, Lemma \ref{lem_osci_est} implies that
         \begin{align*}
             |X_{I_j} |   
             &\leq \beta_{j+1}\sum_{I_{j+1}\subset I_j\cap E_{j+1}}\bigg|\int_{I_{j+1} }e(-\langle \gamma(t),\xi\rangle)dt\bigg| \\
             &\lesssim \beta_{j+1}m_{j+1}\min\{M_{j+1}^{-1},H_{I_{j+1}}^{-1}(\xi)\} \\
             &\lesssim \beta_{j+1} M_j^{-1} 2^{-s_1}.
         \end{align*}
    \end{proof}
    Let $I \subseteq [0,1]$ be an interval such that $I_j \cap I $ is not empty. If we replace $I_j$ in Lemma \ref{lem_XI_s_2=0} by $I_j \cap I$, the proof of Lemma \ref{lem_XI_s_2=0} yields the following corollary, which will be used in the proof of Lemmas \ref{claim_1osc} and \ref{claim_2osc}.
    \begin{corollary}\label{X_I_itrsc}
        For each $\xi \in \widetilde{\Omega}_{j+1}(s_1,s_2)$, We have
        \begin{align}
            |X_{I_j \cap I}|  \lesssim \beta_{j+1} M_j^{-1} 2^{-s_1}.
        \end{align}
    \end{corollary}
        
    Let us turn to the proof of Lemma \ref{claim_1osc}.
    
    Assume that $\xi\not\in  \Omega_{j+1}(s_1,s_2-1)$, then for any $t \in E_{j+1}$, there exists $2 \leq k \leq d$ such that 
    \begin{equation}\label{osc_cond1}
        |\langle \gamma^{(k)} (t), \xi \rangle | \geq M_j 2^{s_1} 2^{(s_2-1)(k-1)}.   
    \end{equation}
    Therefore, for each fixed $\xi \not\in \Omega_{j+1}(s_1,s_2-1)$, we can decompose $[0,1]$ into finitely many closed intervals $\mathcal{J}$ which are disjoint possibly except the endpoints such that for each $I\in \mathcal{J}$, there exists $2\leq k_I\leq d$ which satisfies
    \begin{equation}\label{eq_I_cond_1}
    |\langle \gamma^{(k_I)} (t), \xi \rangle | \geq M_{j+1} 2^{s_1 + (s_2-1)(k_I-1)}, \qquad \forall t \in I    , 
    \end{equation}
    and
    \begin{equation}\label{eq_I_cond_2}
    |\langle \gamma^{(k)} (t), \xi \rangle | \leq M_{j+1} 2^{s_1 + (s_2-1)(k-1)}, \qquad \forall t \in I , \forall 2 \leq k < k_I.    
    \end{equation}
    Since $\langle \gamma(t) ,\xi \rangle$ is a polynomial in $t$ with degree at most $d$, $\#\mathcal{J}$ is bounded uniformly in $j, s_1, s_2$ and $\xi$ and only depends on $d$.
    For each $I \in \mathcal{J}$, let $t_I\in I$ be a number such that $|\langle \gamma^{(1)} (t), \xi \rangle|$ is the minimum on $I$. 
    \begin{lemma}\label{lem_lb_1st_der}
        Assume that $s_2 \neq 0 $, $\xi \in \widetilde{\Omega}_{j+1}(s_1,s_2)$,  and $I \in \mathcal{J}$. Let $m$ be an integer such that $|m| \gtrsim \max\{ M_{j+1}2^{-s_2}, 1\}$. Let $I_{j+1}(t_I+mM_{j+1}^{-1})$ be the interval in $\cI_{j+1}$ which contains $t_I+mM_{j+1}^{-1}$. Then we have
        \[
        |\langle \gamma^{(1)} (t), \xi \rangle| \gtrsim |m| 2^{s_1+s_2}\,,\qquad \forall t\in I_{j+1}(t_I+mM_{j+1}^{-1})\cap I. 
        \]
    \end{lemma}
    \begin{proof}
    For each $t\in I_{j+1}(t_I+mM_{j+1}^{-1})\cap I$, using Taylor expansion up to order $k_I$, we have
    \[
    \begin{split}
    \langle \gamma^{(1)} (t) ,\xi \rangle - \langle \gamma^{(1)} (t_I), \xi \rangle =\sum_{n=2}^{k_I-1}\langle \gamma^{(n)}(t_I),\xi\rangle\frac{(t-t_I)^{n-1}}{(n-1)!}+\langle \gamma^{(k_I)}(t'),\xi\rangle\frac{(t-t_I)^{k_I-1}}{(k_I-1)!},
        \end{split}
    \]
    for some $t'\in I$.
    Since $t_I$ and $t'$ are both in $I$, they satisfy \eqref{eq_I_cond_1} and \eqref{eq_I_cond_2}. Combining with the fact that $|t-t_I|\sim |m|M_{j+1}^{-1}\gtrsim\max \{2^{-s_2},M_{j+1}^{-1}\}$, we obtain that
    \begin{align*}
        |\langle \gamma^{(1)} (t) ,\xi \rangle - \langle \gamma^{(1)} (t_I), \xi \rangle|\sim M_{j+1}2^{s_1} (2^{s_2} mM_{j+1}^{-1})^{k_I-1}.
    \end{align*}
    for sufficiently large $|m| \gtrsim \max\{ M_{j+1}2^{-s_2}, 1 \}$. Since $|\langle \gamma^{(1)} (t_I), \xi \rangle|$ is the minimum and $k_I \geq 2$, the above estimate implies that
    \[
    \begin{split}
    |\langle \gamma^{(1)} (t), \xi \rangle | \gtrsim  |m| 2^{s_1+s_2}.
    \end{split}
    \] 
    Especially, if the signs of $\langle \gamma^{(1)} (t), \xi \rangle$ and $\langle \gamma^{(1)} (t_I), \xi \rangle$ are different, we can use that
    $$
    |\langle \gamma^{(1)} (t),\xi\rangle| \geq \frac{1}{2} |\langle \gamma^{(1)} (t) ,\xi \rangle - \langle \gamma^{(1)} (t_I), \xi \rangle|.
    $$
    \end{proof}  

Now we are ready to prove Lemma \ref{claim_1osc}.
\begin{proof}[Proof of Lemma \ref{claim_1osc}]
    First, let us consider when  $s_2=0$. Lemma \ref{lem_XI_s_2=0} implies that
        \begin{align}
            \sum_{I_j\in\cI_j}|X_{I_j}|^2  &\lesssim \beta_{j+1}^2 M_j^{-2}2^{-2s_1} \#\{I_j \in \cI_j : I_j \subset E_j\} \nonumber \\
            &\lesssim \beta_{j+1}^2\beta_j^{-1} M_j^{-1}2^{-2s_1}  \norm{\mu_j}. \label{sum_XI_s2=0}
        \end{align}
    Therefore, we let $C(E_j,s_1,s_2)$ be a constant multiple of \eqref{sum_XI_s2=0}. Since $\beta_{j+1}^2 \beta_j^{-1} M_j^{-1} \lesssim_\e M_j^{-\alpha+\e}$, Lemma \ref{lem_est_mu_j} implies that a.s. we have \eqref{eq_sum_Xj}, which is
        \begin{align*}
            C(E_j,s_1,s_2) \lesssim_\e M_j^{-\alpha+\e} 2^{-2s_1}\leq M_j^{-\alpha+\e} 2^{-s_1(2-\alpha)}.
        \end{align*}

    Next, we assume that $s_2 \neq 0$. By Cauchy-Schwartz inequality, we have
        \begin{align*}
            \sum_{I_j\in \mathcal{I}_j} |X_{I_j}|^2 =\sum_{I_j\in \mathcal{I}_j} |\sum_{\substack{I\in \mathcal{J}\\I\cap I_j\neq \emptyset}}X_{I_j\cap I}|^2
            \leq \#\mathcal{J}\cdot\bigg(\sum_{I_j\in \mathcal{I}_j}\sum_{\substack{I\in \mathcal{J}\\I\cap I_j\neq \emptyset}} |X_{I_j\cap I}|^2\bigg).
        \end{align*}     
    Since $\# \mathcal{J}$ is bounded, we obtain that
        \begin{align}\label{eq_Xestimate_0}
            \sum_{I_j\in \mathcal{I}_j} |X_{I_j}|^2 \lesssim \sum_{I\in \mathcal{J}}\sum_{\substack{I_j\in \cI_j\\ I\cap I_j\neq \emptyset}}|X_{I_j\cap I}|^2=\sum_{I \in \mathcal{J}}\sum_{|n|\leq M_j} |X_{I_j(t_I+nM_{j}^{-1})\cap I}|^2,
        \end{align}
        where $I_j(t_I+nM_j^{-1})$ is the interval in $\cI_j$ containing $t_I+nM_j^{-1}$.
        
        If $|n| \lesssim \max \{M_j2^{-s_2}, 1\}$, Corollary \ref{X_I_itrsc} implies that
         \begin{align}\label{X_n_small}
             |X_{I_j(t_I + nM_j^{-1}) \cap I} | \lesssim \beta_{j+1} M_j^{-1}2^{-s_1}. 
        \end{align}
         
         If $|n | \gtrsim \max\{ M_j2^{-s_2},1\}$, we use the fact that  
         $$I_j(t_I+nM_j^{-1})\subset \bigcup\limits_{|m-m_{j+1}n|\leq m_{j+1}}I_{j+1}(t_I+mM_{j+1}^{-1}),$$ to get
         \[
        |X_{I_j(t_I + nM_j^{-1})\cap I}| \lesssim \beta_{j+1} \sum_{|m-m_{j+1} n| \leq m_{j+1}}  \left| \int_{I_{j+1}(t_I + mM_{j+1}^{-1})\cap I} e(-\langle \gamma(t), \xi \rangle ) dt \right| .  
        \]
        Since $|m-m_{j+1}n|\leq m_{j+1}$ and $|n| \gtrsim \max\{ M_j 2^{-s_2}, 1\}$ implies that $|m|\gtrsim \max\{M_{j+1}2^{-s_2},1\}$, we apply Lemma \ref{lem_lb_1st_der} to get
        \begin{align*}
            |\langle\gamma^{(1)}(t),\xi\rangle|\gtrsim |m|2^{s_1+s_2}\,,\qquad \forall t\in I_{j+1}(t_I+mM_{j+1}^{-1})\cap I,
        \end{align*}
        which yields that $H_{I_{j+1}(t_I+mM_{j+1}^{-1})\cap I}\gtrsim |m|2^{s_1+s_2}$. By Lemma \ref{lem_osci_est}, we obtain that
        \begin{align*}
            \left| \int_{I_{j+1}(t_I + mM_{j+1}^{-1})\cap I} e(-\langle \gamma(t), \xi \rangle ) dt \right| \lesssim |m|^{-1}2^{-s_1-s_2}.
        \end{align*}
        Hence, we have 
         \begin{equation}\label{eq_X_estimate_3}
         |X_{I_j(t_I + nM_j^{-1})\cap I} | \lesssim \beta_{j+1} |n|^{-1} 2^{-s_1} 2^{-s_2}.    
         \end{equation}

         Now we are ready to establish an upper bound of $\sum_{I_j\in \mathcal{I}_j}|X_{I_j}|^2$. It suffices to sum over the intervals $I_j\subset E_j$ and note that for $r \geq M_j^{-1}$,
         \begin{align*}
            \#\{I_j\subset E_j : I_j\subset B(x,r) \}\lesssim \beta_j^{-1}M_j \mu_j(B(x,r)).
         \end{align*}
         Using \eqref{X_n_small} and \eqref{eq_X_estimate_3}, we obtain that
         \begin{align}
              &\sum_{I\in \mathcal{J}}\sum_{|n|\lesssim \max \{M_j2^{-s_2},1\}} |X_{I_j(t_I+nM_{j}^{-1})\cap I}|^2 \nonumber\\
              & \hspace{2cm}\lesssim  \beta_{j+1}^2 M_j^{-2} 2^{-2s_1}
            \#\{I_j\subset E_j: I_j\subset B(t_I, 2\max \{2^{-s_2},M_j^{-1}\})\}\nonumber\\
            &\hspace{2cm} \leq \beta_{j+1}^2\beta_j^{-1} M_j^{-1} 2^{-2s_1} \mu_j(B(t_I, 2\max \{2^{-s_2},M_j^{-1}\})).\label{sum_C_11}
         \end{align}
         and that
         \begin{align}
            &\sum_{I\in \mathcal{J}}\sum_{|n|\gtrsim \max \{M_j2^{-s_2},1\}} |X_{I_j(t_I+nM_{j}^{-1})\cap I}|^2 \nonumber\\
            &\hspace{2cm} \lesssim \beta_{j+1}^2 2^{-2s_1-2s_2}  \sum_{2^t \gtrsim \max\{ M_j 2^{-s_2},1\} }    2^{-2t}\#\{I_j\subset E_j: I_j\subset B(t_I, 2^{t+1}M_{j}^{-1})\}\nonumber\\
            &\hspace{2cm} \leq \beta_{j+1}^2 \beta_j^{-1} M_j 2^{-2s_1-2s_2}  \sum_{2^t \gtrsim \max\{ M_j 2^{-s_2},1\} }    2^{-2t} \mu_j(B(t_I, 2^{t+1} M_j^{-1}))\label{sum_C_22}
         \end{align}
         By \eqref{eq_Xestimate_0}, we can let $C(E_j,s_1,s_2) $ be a constant multiple of the sum of \eqref{sum_C_11} and \eqref{sum_C_22}.
         
          By Lemma \ref{lem_est_mu_ball_r}, a.s. we obtain that 
         \begin{align*}
            \beta_{j+1}^2\beta_j^{-1} M_j^{-1} 2^{-2s_1} \mu_j(B(t_I, 2\max \{2^{-s_2},M_j^{-1}\})) \lesssim_\e M_j^{-\alpha+\e}2^{-2s_1} (\max\{2^{-s_2}, M_j^{-1}\})^\alpha
         \end{align*}
         and that
        \begin{align*}
            \beta_{j+1}^2 \beta_j^{-1} M_j 2^{-2s_1-2s_2}  &\sum_{2^t \gtrsim \max\{ M_j 2^{-s_2},1\} }    2^{-2t} \mu_j(B(t_I, 2^{t+1} M_j^{-1}))\\
            &\hspace{2cm}\lesssim_\e M_j^{2-2\alpha+\e} 2^{-2s_1-2s_2} \sum_{2^t \gtrsim \max\{ M_j 2^{-s_2},1\} }   2^{ -(2-\alpha)t}\\
            &\hspace{2cm}\lesssim_\e M_j^{2-2\alpha+\e} 2^{-2s_1-2s_2} (\max\{ M_j 2^{-s_2},1\} )^{ -2+\alpha}\\
            &\hspace{2cm}\lesssim_\e M_j^{-\alpha+\e} 2^{-2s_1-2s_2} (\max\{  2^{-s_2},M_j^{-1}\} )^{ -2+\alpha}.
        \end{align*}
        Combining two estimates and the fact that $ 2^{s_2}\leq 2^{s_1}M_j$, a.s. we get 
        \begin{align*}
            C(E_j,s_1,s_2) &\lesssim_\e M_j^{-\alpha+\e} 2^{-2s_1} \left( \max \{2^{-s_2\alpha} , M_j^{-\alpha}\}+ 2^{-2s_2} (\max\{ 2^{-s_2},M_j^{-1}\} )^{ -2+\alpha}\right) \\
            &\lesssim_\e M_j^{-\alpha+\e} 2^{-s_1(2-\alpha)} 2^{-s_2\alpha}.
        \end{align*}
\end{proof}

In order to prove Lemma \ref{claim_2osc} We need to divide $[0,1]$ into sets where Lemma \ref{lem_lb_1st_der} is applicable and where it is not. For fixed $\xi\in \widetilde{\Omega}_{3,j+1}(s_1,s_2)$, we define
    \begin{align*}
        \Lambda_1&=\big\{ t\in [0,1]: \exists k\geq 2 \text{ such that }
|\langle \gamma^{(k)}(t) , \xi \rangle | \geq 10^{-8}(M_{j} 2^{s_1}) 2^{s_2 (k-1)}\big\}
    \end{align*}
    and $\Lambda_2 := [0,1] \backslash \Lambda_1$. 

Assume that $\xi \in \widetilde{\Omega}_{3,j+1}(s_1,s_2)$ for $s_2 \neq 0$. Observe that since $\xi \in \Omega_{3,j+1}(s_1,s_2)$, there exists $t_0$ such that
    \begin{align}\label{cond_t_00}
        |\langle \gamma^{(k)}(t_0), \xi \rangle | \leq M_{j+1} 2^{s_1} 2^{s_2(k-1)},\qquad \forall 1 \leq k \leq 3.
    \end{align}
    Also, since $\xi \notin \Omega_{3,j+1}(s_1,s_2-1)$, we have 
    \begin{equation}\label{cond_t_0}
      |\langle \gamma^{(2)}(t_0), \xi \rangle | > M_{j+1} 2^{s_1} 2^{(s_2-1)} \qquad \text{or} \qquad  |\langle \gamma^{(3)}(t_0), \xi \rangle | > M_{j+1} 2^{s_1} 2^{2(s_2-1)}.  
    \end{equation}
    Note that $ |\langle \gamma^{(1)}(t_0), \xi \rangle | > M_{j+1} 2^{s_1} $ cannot hold because of \eqref{cond_t_00}. Hence, $t_0 \in  \Lambda_1$ and $\Lambda_1$ is nonempty. Since cannot use Lemma \ref{lem_lb_1st_der} on $\Lambda_2$, we need an additional lemma.
    \begin{lemma}\label{I_2_est}
        For $\xi \in \widetilde{\Omega}_{3,j+1}(s_1,s_2)$ with $s_2 \neq 0$, assume that $$M_{j+1}2^{s_1+s_2+s_3-1} \leq |\xi_3|\leq  M_{j+1} 2^{s_1+s_2+s_3}$$ 
        where $ 2^{s_3} \leq 10^{-8} 2^{s_2}$. If $2^{s_3} < 100^{-1}$, then $\Lambda_2$ is empty. If $100^{-1} \leq 2^{s_3}$, then for $t_0 \in \Lambda_1$ defined above and for any $t \in \Lambda_2$, we have
        $$
        |t-t_0| \sim 2^{-s_3}
        $$
        and
        \begin{equation}\label{gamma1_lb_d=3}
         |\langle \gamma^{(1)} (t), \xi \rangle | \gtrsim M_{j+1}2^{s_1+s_2-s_3}.  
        \end{equation}
    \end{lemma}
    \begin{proof}
    Since $\xi_3 \neq 0$, $\langle \gamma(t), \xi\rangle$ is a polynomial in $t$ with degree $3$ and it implies that
    \begin{equation}\label{est_xi3}
      |\langle \gamma^{(3)} (t), \xi \rangle | =6 |\xi_3| \sim M_{j+1}2^{s_1+s_2+s_3} \qquad \forall t \in [0,1].  
    \end{equation}
    Since $2^{s_3} \leq 10^{-8} 2^{s_2}$, \eqref{cond_t_0} implies that
    \begin{equation}\label{est_t_0}
      M_{j+1} 2^{s_1} 2^{(s_2-1)} < |\langle \gamma^{(2)}(t_0), \xi \rangle | \leq M_{j+1} 2^{s_1} 2^{s_2}.  
    \end{equation}
    Also, since $t  \in \Lambda_2$,  we have
    \begin{equation}\label{est_t}
        |\langle \gamma^{(2)} (t), \xi \rangle| \leq 10^{-8}M_{j+1} 2^{s_1+s_2}.
    \end{equation}
    
    For each $t$, observe that
    \[
    |\langle \gamma^{(2)} (t), \xi \rangle - \langle \gamma^{(2)} (t_0), \xi \rangle | = |\langle \gamma^{(3)} (t_0), \xi \rangle  (t-t_0)| =6 |\xi_3| |t-t_0|.
    \]
    Combining \eqref{est_xi3}, \eqref{est_t_0} and \eqref{est_t}, we obtain that
    \begin{align}\label{range_t}
        100^{-1} 2^{-s_3} \leq |t-t_0 | \leq  100\cdot  2^{-s_3}.
    \end{align}
    Since $\Lambda_2 \subseteq [0,1]$, $\Lambda_2$ is empty if $2^{s_3} < 100^{-1}$. If $100^{-1} \leq 2^{s_3}$, then we have $|t-t_0| \sim 2^{-s_3}$.
    
    For the lower bound of $|\langle \gamma^{(1)} (t), \xi \rangle |$, we use Taylor expansion \eqref{eq_taylor} and obtain that
    \[
    \langle \gamma^{(1) } (t_0) , \xi \rangle =  \langle \gamma^{(1)} (t),\xi \rangle + \frac{\langle \gamma^{(2)} (t),\xi \rangle}{2}  (t_0-t) + \frac{\langle \gamma^{(3)} (t),\xi \rangle}{6} (t_0-t)^2.
    \]
    Using \eqref{est_xi3}, \eqref{est_t} and \eqref{range_t}, we get
    \[
    \begin{split}
     |\langle \gamma^{(1)} (t_0 ), \xi \rangle| &\leq M_{j+1}2^{s_1}\\
     |\langle \gamma^{(2)} (t), \xi \rangle (t_0-t)| &\leq 10^{-6} M_{j+1} 2^{s_1+s_2-s_3}  \\
    \end{split}
    \]
    while
    \[
    |\langle \gamma^{(3)} (t), \xi \rangle (t_0-t)^2| \geq 10^{-4}M_{j+1} 2^{s_1+s_2-s_3}
    \]
    Thus, we obtain \eqref{gamma1_lb_d=3} as desired. 
    
    \end{proof}
    \begin{remark}
        The proof of Lemma \ref{I_2_est} works only when $d=3$, since we used that $|\langle \gamma^{(3)} (t), \xi \rangle |$ is a constant.
    \end{remark}
    
    Finally, we prove Lemma \ref{claim_2osc}.
    \begin{proof}[Proof of Lemma \ref{claim_2osc}]
    If $s_2=0$, the proof is the same as Lemma \ref{claim_1osc}.
    
    Thus, it suffices to assume that $s_2 \neq 0$. First, we divide $\Lambda_1$ into finitely many intervals $I$ such that there exists $k_I= 2$ or $3$ which satisfies 
        \begin{equation}\label{eq_I_cond_1_d=3}
    |\langle \gamma^{(k_I)} (t), \xi \rangle | \geq 10^{-8} M_j 2^{s_1 + (s_2-1)(k_I-1)}, \qquad \forall t \in I    , 
    \end{equation}
    and
    \begin{equation}\label{eq_I_cond_2_d=3}
    |\langle \gamma^{(k)} (t), \xi \rangle | \leq 10^{-8}   M_j 2^{s_1 + (s_2-1)(k-1)}, \qquad \forall t \in I , \forall 2 \leq k < k_I.    
    \end{equation}
    We denote the set of such intervals by $\mathcal{J}_1$. Note that $\# \mathcal{J}_1$ is uniformly bounded in $\xi$, $s_1$ and $s_2$. The set $\Lambda_2$ is a connected interval since $\langle \gamma(t),\xi \rangle $ is a polynomial in $t$ with degree at most $3$. Thus, we do not need to divide $\Lambda_2$ further. Since $\# \mathcal{J}_1 \lesssim 1$, we have
    \[
    \begin{split}
        \sum_{I_j\in \mathcal{I}_j} |X_{I_j}|^2 &=\sum_{I_j\in \mathcal{I}_j} |\sum_{\substack{I\in \mathcal{J}_1\\I\cap I_j\neq \emptyset}}X_{I_j\cap I}|^2 + \sum_{I_j\in \mathcal{I}_j} |X_{I_j\cap \Lambda_2}|^2.
    \end{split}
    \]
    For the first sum, we can use the same argument from Lemma \ref{claim_1osc} with \eqref{eq_I_cond_1} and \eqref{eq_I_cond_2} replaced by \eqref{eq_I_cond_1_d=3} and \eqref{eq_I_cond_2_d=3} respectively. Thus, there exists $C_1(E_j,s_1,s_2)$ such that
    $$
    \sum_{I_j\in \mathcal{I}_j} |\sum_{\substack{I\in \mathcal{J}_1\\I\cap I_j\neq \emptyset}}X_{I_j\cap I}|^2 \leq C_1(E_j,s_1,s_2)
    $$
    conditioned on $E_j$ and a.s. we have
    $$
    C_1(E_j,s_1,s_2) \lesssim_\e M_j^{-\alpha+\e} 2^{-s_1(2-\alpha)} 2^{-s_2\alpha}
    $$
    uniformly in $\xi$ and $E_j$.
    
    For the given  $\xi \in \widetilde{\Omega}_{3,j+1}(s_1,s_2)$, let us assume that $|\xi_3| \sim M_{j+1} 2^{s_1+s_2+s_3}$ where $2^{s_3} \leq 2^{s_2} $. We will find an estimate on the second sum which is uniform on $s_3$. If $s_3 < 100^{-1} $, $\Lambda_2$ is empty by Lemma \ref{I_2_est}. If $s_3 \sim s_2$, then 
    \begin{equation*}
    |\langle \gamma^{(3)} (t), \xi \rangle | \sim M_j 2^{s_1 + 2s_2}, \qquad \forall t \in [0,1]  , 
    \end{equation*}
    and since $t \in \Lambda_2$, we have
    \begin{equation*}
    |\langle \gamma^{(2)} (t), \xi \rangle | \lesssim M_j 2^{s_1 + s_2}, \qquad \forall t \in [0,1]. 
    \end{equation*}
    They corresponds to \eqref{eq_I_cond_1} and \eqref{eq_I_cond_2} respectively with $k_I=3$ and $I=[0,1]$. Thus, we can use the argument from Lemma \ref{claim_1osc} again. Therefore, it suffices to consider when $100^{-1} \leq 2^{s_3} \leq 10^{-8} 2^{s_2}$.

    If $I_j(t_0+nM_j^{-1}) \cap \Lambda_2 \neq \emptyset$ for $|n| \gtrsim 1$, Lemmas \ref{lem_osci_est} and  \ref{I_2_est} imply that $ |n| \sim  M_j2^{-s_3}$ and
        \[
        |X_{I_j(t_0+nM_j^{-1}\cap \Lambda_2)}| \lesssim \beta_{j+1} M_j^{-1} 2^{-s_1-s_2+s_3}.
        \]
    Since $|n|M_j^{-1} \sim 2^{-s_3} \gtrsim M_j^{-1}$, the number of $I_j$ such that $I_j(t_0+nM_j^{-1}) \cap \Lambda_2 \neq \emptyset$ is bounded by $\beta_j^{-1} M_j \mu_j( B(t_0, C2^{-s_3})) $ for sufficiently large $C$, say $200$.

    If $|n| \lesssim 1 $, we use Corollary \ref{X_I_itrsc} and get 
    \[
        |X_{I_j(t_0+nM_j^{-1}\cap \Lambda_2)}| \lesssim \beta_{j+1} M_j^{-1} 2^{-s_1}.
    \]
    Therefore,
    \begin{align*}
        \sum_{I_j\in \cI_j} |X_{I_j\cap \Lambda_2}|^2 \lesssim  &\beta_{j+1}^2 M_j^{-2} 2^{-2s_1}\\
        &+  \beta_{j+1}^2\beta_j^{-1} M_j^{-1} 2^{-2s_1-2s_2+2s_3} \mu_j( B(t_0, C2^{-s_3}))w_{j}(s_3)
    \end{align*}
    where $w_j(s_3) = 1 $ if $2^{s_3} \leq C^{-1} M_j$ and $w_j(s_3) =0$ otherwise. 
    Thus, we let $C_2(E_j,s_1,s_2)$ be a constant multiple of the right-hand side.

    Now, we find an upper bound of $C_2(E_j,s_1,s_2)$. Since $2^{s_2} \leq M_j 2^{s_1}$, we get
    \begin{align*}
        \beta_{j+1}^2 M_j^{-2} 2^{-2s_1} \lesssim_\e M_j^{-2\alpha+\e}2^{-2s_1}\leq M_j^{-\alpha+\e}2^{-2s_1+(s_1-s_2)\alpha}.
    \end{align*}
    Using Lemma \ref{lem_est_mu_ball_r} and the fact that $s_3 \leq s_2$, a.s. we have
    \begin{align*}
        \beta_{j+1}^2\beta_j^{-1} M_j^{-1} 2^{-2s_1-2s_2+2s_3} \mu_j( B(t_0, C2^{-s_3}))w_j(s_3) 
        &\lesssim_\e M_j^{-\alpha +\e} 2^{-2s_1-2s_2+(2-\alpha)s_3}\\
        &\lesssim_\e M_j^{-\alpha +\e} 2^{-2s_1- s_2\alpha}
    \end{align*}
    Combining two estimates, a.s. we obtain that 
    \begin{align*}
            C_2(E_j,s_1,s_2)  \lesssim_\e M_j^{-\alpha+\e} 2^{-s_1(2-\alpha)} 2^{-s_2\alpha}.
        \end{align*}
    Let $C(E_j,s_1,s_2) := C_1(E_j,s_1,s_2)+C_2(E_j,s_1,s_2)$. Then, $C(E_j,s_1,s_2)$ satisfies the desired properties. 
    \end{proof}

\bibliographystyle{abbrv}
\bibliography{ref} 

\end{document}